\documentclass[12pt]{article}
\usepackage{amsmath,amssymb,amsthm,mathtools}
\usepackage{xcolor}
\allowdisplaybreaks

\usepackage{booktabs,array,fullpage}
\usepackage{enumitem}
\usepackage{microtype}
\usepackage{hyperref,todonotes}

\hypersetup{
  colorlinks=true,
  linkcolor=blue,
  citecolor=blue,
  urlcolor=blue,
  pdftitle={From Gaussian to stable and back: Non-Pareto points in growing dimension}
}

\newtheorem{theorem}{Theorem}[section]
\newtheorem{proposition}[theorem]{Proposition}
\newtheorem{lemma}[theorem]{Lemma}
\newtheorem{corollary}[theorem]{Corollary}

\theoremstyle{remark}
\newtheorem{remark}[theorem]{Remark}

\numberwithin{equation}{section}

\newcommand{\E}{\mathbb E}
\newcommand{\R}{\mathbb R}
\renewcommand{\P}{\mathbf P}
\newcommand{\Prob}[1]{\mathbf{P}\left\{#1\right\}}
\newcommand{\sI}{\mathcal{I}}

\newcommand{\sF}{\mathcal{F}}

\newcommand{\sO}{\mathcal{O}}
\newcommand{\Var}{\operatorname{Var}}
\newcommand{\Cov}{\operatorname{Cov}}
\newcommand{\card}{\operatorname{card}}

\newcommand{\one}{\mathbf 1}
\newcommand{\dd}{\,\mathrm d}
\newcommand{\eps}{\varepsilon}

\newcommand{\wto}{\xrightarrow{d}}
\newcommand{\pto}{\xrightarrow{\P}}

\definecolor{revisiongreen}{RGB}{0,110,45}

\allowdisplaybreaks[4]
\usepackage{version}

\title{From Gaussian to stable and back:\\Non-Pareto points in growing dimension}
\author{Chinmoy Bhattacharjee
  \footnote{Department of Mathematics, Universit\"at Hamburg,
    Bundesstra{\ss}e 55, 20146 Hamburg, Germany, chinmoy.bhattacharjee@uni-hamburg.de},
  Andrii Ilienko
  \footnote{Faculty of Pure and Applied Mathematics, Wrocław
    University of Science and Technology, 50-370 Wroc{\l}aw, Poland,
    andrii.ilienko@pwr.edu.pl},
  Bochen Jin
  \footnote{Institute of Mathematical Statistics and Actuarial
    Science, University of Bern, Alpeneggstr. 22, CH-3012, Bern,
    Switzerland, bochen.jin@unibe.ch},
  Ilya Molchanov
  \footnote{Institute of Mathematical Statistics and Actuarial
    Science, University of Bern, Alpeneggstr. 22, CH-3012, Bern,
    Switzerland, ilya.molchanov@unibe.ch}
}
\date{\today}

\begin{document}

\maketitle

\begin{abstract}
  We study the number of non-Pareto points in two high-dimensional models: a homogeneous Poisson point process and a fixed-size sample of independent uniform points in the unit cube. A point is non-Pareto if some other sample point is coordinatewise less than or equal to it. As the dimension grows on the logarithmic scale of the intensity or sample size, the asymptotic behavior exhibits a surprising sequence of transitions: from Gaussian to stable fluctuations, back to Gaussian, then to a Poisson limit at the sparsity threshold, followed by the eventual disappearance of non-Pareto points.
  
  We identify the critical thresholds, normalizations, and limiting laws in both models. Notably, the return to Gaussianity occurs while the stability index remains strictly below $2$ and is driven by a change in the dominant fluctuation mechanism. Moreover, in part of the returning Gaussian regime, the natural normalization is the square root of the expected count rather than its standard deviation. At the stable--Gaussian transition, the two mechanisms can coexist, giving rise to a sum of independent stable and Gaussian random variables.
\end{abstract}




\section{Introduction}

A point of a finite subset of $\R^d$ is called Pareto-minimal (or simply Pareto) if no other point of the set is coordinatewise less than or equal to it. The set of Pareto points, called the Pareto front, plays a central role in multicriteria optimization~\cite{Deb2009,Ehrgott2005}. Related coordinatewise comparisons arise in a wide range of other settings, including maximal-vector algorithms and their applications to convex hulls~\cite{BentleyKungSchkolnickThompson1978,D80}, skyline queries in databases~\cite{BorzsonyiKossmannStocker2001}, etc. As the dimension grows, coordinatewise comparability between random points becomes increasingly rare. This motivates studying the complement of the Pareto front --- the non-Pareto points.

We consider two models of random points in $Q_d=[0,1]^d$: a
homogeneous Poisson point process $\eta_s$ of intensity $s$ and a
sample of $n$ independent uniform points. We allow the dimension $d$
to depend on $s$ or $n$ and denote the corresponding numbers of
non-Pareto points by $N_s$ and $N_n^{\rm bin}$. For two independent
uniform points $X,Y\in Q_d$, we have $\Prob{X\le
  Y}=2^{-d}$. Consequently, the expected numbers of ordered pairs of
distinct points $(x,y)$ satisfying $x\le y$ are $s^2 2^{-d}$ and
$n(n-1)2^{-d}$ in the Poisson and binomial models, respectively. This
competition between the growing number of potential pairs and the
exponentially decreasing probability of coordinatewise comparability
naturally singles out the logarithmic scale for the dimension: $d/\log
s\to\beta$ and $d/\log n\to\beta$ in the respective models.
We also note that the homogeneous intensity measure in the Poisson
model can be replaced by any product measure with continuous marginals
and total mass $s$, while the uniform sampling distribution in the
binomial model by any distribution with independent continuous
marginals. Since coordinatewise quantile transforms preserve the
partial order, all our results remain valid in this general setting.

The principal finding of this paper is a surprising non-monotone
sequence of limiting laws. As $\beta$ increases, Gaussian fluctuations
give way to stable ones and then reappear, followed by a Poisson limit
at the sparsity threshold and the eventual disappearance of non-Pareto
points. Moreover, the return to Gaussianity is not due to the
stability index approaching $2$, but to a change in the dominant
fluctuation mechanism. This change also leads to an unusual
normalization: in the lower part of the returning Gaussian regime, the
centered count of non-Pareto points has a Gaussian limit when
normalized by the square root of its expectation rather than by its
standard deviation. 

The study of random Pareto points in fixed dimension goes back to
Barndorff-Nielsen and
Sobel~\cite{BarndorffNielsenSobel1966}. Subsequent results include
variance asymptotics~\cite{BaiChaoHwangLiang1998} and Gaussian limit
theorems for the number of Pareto points, both in fixed-size
samples~\cite{BDHT05,BaiHwangLiangTsai2001,BarbourXia2001} and in
Poisson samples~\cite{B00}. More recently, Bhattacharjee and
Molchanov~\cite{BM22} obtained quantitative Gaussian approximations in
the Poisson cube model by extending stabilization
methods~\cite{MR3910021,MR2201885}. In growing dimension, Jacobovic
and Zuk~\cite{JZ23} investigated a related but different question: the
transition in the probability that a given sample point is Pareto. 

The existence of non-Pareto points is closely related to chains in the
coordinatewise partial order: a non-Pareto point exists if and only if
the sample contains a chain of length two. Bollob\'as and
Brightwell~\cite{BollobasBrightwell1995} established sharp thresholds
for the existence of chains of fixed length in a sample of $n$
independent uniform points. In particular, the sparsity threshold for
two-point chains is $d_n=\beta_{\mathrm{sp}}\log n$, where
$\beta_{\mathrm{sp}}=2/\log 2\approx 2.885$. More recently, Ilienko
and Jin~\cite{IlienkoJin2026} obtained Poisson point-process limits
for finite-dimensional coordinate projections of non-Pareto points in
the critical regime $d_n-\beta_{\mathrm{sp}}\log n\to c\in\R$. The
limiting processes have total intensity $2^{-c}$, yielding in
particular a Poisson limit with mean $2^{-c}$ for the number of
non-Pareto points. The same limits hold in the homogeneous Poisson
model, with $n$ replaced by the intensity $s$. These results describe
the transition to sparsity but leave open the fluctuation behavior
below the threshold. 

Gaussian--stable transitions arising from light-tailed input have
been studied extensively for sums of random exponentials.
Ben Arous, Bogachev and Molchanov~\cite{BenArousBogachevMolchanov2005}
identify growth regimes in which such sums satisfy laws of large
numbers, Gaussian limit theorems, or stable limit theorems.
Closely related results for sums of independent random products
were obtained by Cranston and Molchanov~\cite{CranstonMolchanov2005}.
Kabluchko~\cite{Kabluchko2011GeometricLevy} established functional
limit theorems for sums of independent geometric L\'evy processes,
including critical Gaussian limits whose variances depend on
second-order growth conditions. These results demonstrate how
rare large summands can generate stable fluctuations even when
the underlying variables possess finite exponential moments.

We now focus on the Poisson model, for which our main results provide
a detailed distributional phase diagram below the sparsity threshold.
The binomial model exhibits the same principal sequence of limiting
laws, with the corresponding results obtained by de-Poissonization in
Section~\ref{sec:de-poiss-binom}. 

We identify two additional critical values, $\beta_L\approx 2.315$ and
$\beta_c\approx 2.462$, defined explicitly in
Section~\ref{sec:model-main-stat}. For $\beta\in[0,\beta_L)$, the
centered count $N_s-\E N_s$ has a Gaussian limit after normalization
by its standard deviation. In the intermediate range
$\beta\in(\beta_L,\beta_c)$, a smaller normalization yields a centered
spectrally positive stable law with index $\alpha_\beta$. Remarkably,
as $\beta$ increases across this interval, $\alpha_\beta$ decreases
from $2$ to $2/\log 2-1\approx1.885$, so the stable phase covers only
a narrow range of indices close to $2$. For
$\beta\in(\beta_c,\beta_{\mathrm{sp}})$, Gaussian fluctuations
reappear, now after normalization by $\sqrt{\E N_s}$. This
normalization is asymptotically equivalent to the standard deviation
for $\beta>\beta_1\approx2.466$; in the narrow interval
$(\beta_c,\beta_1)$, normalization by the standard deviation instead
yields a degenerate weak limit. In particular, even the second moments
fail to converge to that of the limit.

The critical boundaries exhibit further unusual behavior. At
$\beta_L$, normalization by the standard deviation yields Gaussian
limits whose variances can take any value in $[0,1]$, depending on the
second-order behavior of the dimension. At $\beta_c$, the stable and
Gaussian fluctuation mechanisms can coexist, yielding a limiting
distribution given by the sum of independent stable and Gaussian
random variables.

The different limiting laws arise from a competition between distinct
sources of fluctuations. A predictable decomposition separates the
count into a counting martingale and a random compensator, whose
Wiener--It\^o chaos expansion isolates its first chaos. Even within
the lower Gaussian regime, the dominant component changes at an
additional threshold $\beta_0\approx0.305$: below $\beta_0$, the
fluctuations are governed by the total number of points in the Poisson
process, whereas for $\beta\in(\beta_0,\beta_L)$, the first chaos of
the compensator dominates the remaining terms and has a Gaussian
limit. For $\beta\in(\beta_L,\beta_c)$, the same first chaos continues
to determine the fluctuations of the full count, but now on a smaller
scale, with rare large contributions generating a stable limit. In
fact, the first chaos retains its stable limit throughout the much
wider interval $(\beta_L,\beta_{\mathrm{sp}})$, with the stability
index continuing to decrease as $\beta$ increases. Beyond $\beta_c$,
however, the counting martingale becomes dominant in the fluctuations
of the full count, restoring the Gaussian limit. The abrupt return to
Gaussianity at $\beta_c$ is therefore caused by a change in the
dominant component, even though the stability index remains strictly
below $2$. Thus the two Gaussian regimes arise from fundamentally
different probabilistic mechanisms. We also show that the upper
Gaussian regime extends to the dense side of the sparsity threshold,
that is, when $d_s/\log s\to\beta_{\mathrm{sp}}$ and
$d_s-\beta_{\mathrm{sp}}\log s\to-\infty$.

The compensator first chaos is also where the connection with the
existing literature on sums of random exponentials becomes explicit.
Up to centering and a deterministic factor, its product-kernel
approximation, introduced in Subsection~\ref{sec:zero-saddle}, is a
sum of $\operatorname{Poisson}(s)$ many independent products of $d_s$
i.i.d.\ standard uniform random variables, that is, a Poissonized
version of the sums of random products studied in
\cite{CranstonMolchanov2005,Kabluchko2009RandomProducts}. For the
uniform-product example in
\cite[Section~1.5]{Kabluchko2009RandomProducts}, the Gaussian--stable
threshold for the exponential growth rate $c$ of the number of
summands relative to the number of factors is $c_2=\log3-2/3$. In our
setting, we have $c = 1/\beta$ and our threshold $\beta_L$ satisfies
$\beta_L=c_2^{-1}$. The corresponding stable-index equation yields the
same index $\alpha_\beta$ as in the present paper. One of our
contributions lies in connecting this independent-product mechanism to
the dependent geometric count of non-Pareto points. This requires
controlling the product-kernel approximation and the remaining
fluctuation terms on the relevant scales. We identify a further
transition at $\beta=\beta_c$, caused by the competition between the
compensator first chaos and the counting martingale; this return back
to the Gaussian behaviour has no counterpart for sums of independent
products.

The proofs combine precise asymptotic estimates with limit theorems
for Poisson integrals. High-dimensional integral estimates yield the
required asymptotics for the mean and variance contributions, as well
as bounds on the higher-order chaos remainder, while Gamma lower-tail
asymptotics determine the stable normalization and stability index. In
the lower Gaussian regime, we establish asymptotic normality using
Lyapunov-type moment estimates and a Lindeberg-type analysis of large
jumps. For the upper Gaussian regime and the transition at $\beta_c$,
we prove a general joint limit theorem for predictable and
deterministic Poisson integrals. This result yields both the returning
Gaussian limit and, at the transition, the asymptotic independence of
the Gaussian and stable components.

The rest of the paper is organized as
follows. Section~\ref{sec:model-main-stat} introduces the model and
presents the main results. Section~\ref{sec:compensator-chaos}
develops the martingale and Poisson chaos decompositions and derives
an exact variance identity. Section~\ref{sec:prod-minim-lapl}
establishes asymptotics for the relevant high-dimensional integrals,
which are then used in Section~\ref{sec:comp-first-chaos} to compare
the variance contributions and control the higher-order chaos
remainder. Section~\ref{sec:gauss-rang-beta_0b} proves the lower
Gaussian limit and examines the boundary at $\beta_L$. The stable
regime is treated in Section~\ref{sec:above-lind-trans}, and the upper
Gaussian regime, including its extension to the dense side of the
sparsity threshold, in
Section~\ref{sec:above-betac}. Section~\ref{sec:bound-behav-at-1}
analyzes the stable--Gaussian transition at $\beta_c$, including the
possible mixed limiting laws. Finally,
Section~\ref{sec:de-poiss-binom} establishes the corresponding limit
theorems for the binomial model by de-Poissonization.

\section{The Poisson model: setup and main results}
\label{sec:model-main-stat}

Let $d_s\in\mathbb N$ be a dimension sequence satisfying
\begin{displaymath}
	\frac{d_s}{\log s}\to\beta\in[0,\infty)
	\qquad\text{as }s\to\infty.
\end{displaymath}
Unless otherwise stated, all subsequent limits are likewise taken as
$s\to\infty$. For each $s$, let $\eta_s$ be a homogeneous Poisson
point process on $Q_{d_s}=[0,1]^{d_s}$ with intensity $s$. The
processes $\eta_s$ for different values of $s$ need not be coupled.

For $x,y\in\R^{d_s}$, we write $x\le y$ if $x_i\le y_i$ for all $i=1,\ldots,d_s$. We also use the notation
\begin{displaymath}
  |x|=\prod_{i=1}^{d_s}x_i,\qquad
  [x,y]=[x_1,y_1]\times\cdots\times[x_{d_s},y_{d_s}],
  \qquad
  [x,y)=[x,y]\setminus\{y\}.
\end{displaymath}
A point $x\in\eta_s$ is called Pareto-minimal (or simply Pareto) if no other point $y$ of $\eta_s$ satisfies $y\le x$, that is, if $\eta_s([0,x))=0$. Thus the number of non-Pareto points is given by
\begin{displaymath}
  N_s=\card\{x\in\eta_s:x\text{ is not Pareto}\}
  =\sum_{x\in\eta_s}\one_{\{\eta_s([0,x))\geq1\}}.
\end{displaymath}

We now present our main results, which involve the following critical thresholds:
\begin{align*}
  \beta_L
  &=\frac{3}{3\log3-2}\approx2.315,
  &\qquad \beta_c
  &=\left(\frac{\log2}{2}-1-\log\frac{\log2}{2}\right)^{-1}\approx2.462,\\
  \beta_1
  &=\frac1{\log(3/2)}\approx2.466,
  &\qquad \beta_{\mathrm{sp}}
  &=\frac2{\log2}\approx2.885.
\end{align*}

We begin with the lower Gaussian regime.

\begin{theorem}
  \label{thm:main-gaussian}
  If $d_s/\log s\to\beta\in[0,\beta_L)$, then
  \begin{displaymath}
    \frac{N_s-\E N_s}{\sqrt{\Var N_s}}
    \wto N(0,1).
  \end{displaymath}
\end{theorem}

At the boundary $\beta_L$, the first-order condition
$d_s/\log s\to\beta_L$ alone does not determine a unique
limiting distribution. The limit depends on the deviation
$d_s-\beta_L\log s$ on the scale
$\sqrt{d_s}\asymp\sqrt{\log s}$.
Define
\begin{equation}
  \label{eq:vartheta-s}
  \vartheta_s=
  \frac{3(d_s-\beta_L\log s)}
  {2\beta_L\sqrt{d_s}}
\end{equation}
and let $\Phi$ denote the standard normal distribution function.

\begin{theorem}
  \label{thm:limtbetaL}
  If $d_s/\log s\to\beta_L$ and
  $\vartheta_s\to\vartheta\in[-\infty,+\infty]$, then
  \begin{equation*}
    \frac{N_s-\E N_s}{\sqrt{\Var N_s}}
    \wto N\bigl(0,\Phi(-\vartheta)\bigr).
  \end{equation*}
\end{theorem}

We next turn to the stable regime.  Let $\sI(u)=u-1-\log u$, $u>0$, be
the Cram\'er rate function of the $\operatorname{Exp}(1)$
distribution. For $\beta>0$, let $r_\beta\in(0,1)$ be the unique
solution of $\sI(r_\beta)=1/\beta$, and set
$\alpha_\beta=1/r_\beta-1$.  In particular, as shown at the beginning
of Subsection~\ref{sec:gamma-normalization}, $\alpha_\beta\in(1,2)$
for $\beta\in(\beta_L,\beta_{\mathrm{sp}})$.

Let $\Gamma_{d_s}\sim\operatorname{Gamma}(d_s,1)$.  For $s>1$, let
$t_s^\star$ be the unique solution of
$\P\{\Gamma_{d_s}\le t_s^\star\}=1/s$, and define
\begin{equation}\label{eq:bs}
  b_s=s e^{-t_s^\star}.
\end{equation}

\begin{theorem}
  \label{thm:stable-full}
  If $d_s/\log s\to\beta\in(\beta_L,\beta_c)$,
  then
  \begin{displaymath}
    \frac{N_s-\E N_s}{b_s}
    \wto Z_{\alpha_\beta},
  \end{displaymath}
  where $Z_{\alpha_\beta}$ is the centered
  spectrally positive $\alpha_\beta$-stable
  random variable with characteristic exponent
  \begin{displaymath}
    \log\E e^{itZ_{\alpha_\beta}}
    =\int_0^\infty
    (e^{itu}-1-itu)
    \alpha_\beta u^{-1-\alpha_\beta}\dd u,
    \qquad t\in\R.
  \end{displaymath}
\end{theorem}

At the critical value $\beta_c$, the stable and
Gaussian components can contribute on comparable
scales. The resulting limit depends on a finer
centering of the dimension. Write
$\alpha_c=\alpha_{\beta_c}=2/\log2-1\approx1.885$
and define
\begin{equation}\label{eq:dsstar}
  d_s^\star=\beta_c\log s
  -\frac{\beta_c}{2}\log\log s
  -\beta_c\log\bigl(
  (1-r_{\beta_c})\sqrt{2\pi\beta_c}\bigr).
\end{equation}

\begin{theorem}
  \label{thm:boundary-trichotomy}
  Suppose that $d_s/\log s\to\beta_c$.
  \begin{enumerate}[label=\textup{(\roman*)}]
  \item \label{it:i} If $d_s-d_s^\star\to-\infty$, then
    \begin{displaymath}
      \frac{N_s-\E N_s}{b_s}
      \wto Z_{\alpha_c}.
    \end{displaymath}
  \item \label{it:ii} If $d_{s_n}-d_{s_n}^\star\to x\in\R$
  along some sequence $s_n\to\infty$, then
  \begin{displaymath}
  	\frac{N_{s_n}-\E N_{s_n}}{b_{s_n}}
  	\wto Z_{\alpha_c}
  	+e^{x/(\alpha_c\beta_c)}G,
  \end{displaymath}
  where $G\sim N(0,1)$ is independent of $Z_{\alpha_c}$.
  \item If $d_s-d_s^\star\to+\infty$, then
    \begin{displaymath}
      \frac{N_s-\E N_s}{\sqrt{\E N_s}}
      \wto N(0,1).
    \end{displaymath}
  \end{enumerate}
\end{theorem}

The restriction to a subsequence in \ref{it:ii}
is necessary, since $d_s$ is integer-valued
while $d_s^\star$ varies continuously.

Beyond $\beta_c$, Gaussian fluctuations reappear throughout the range
below the sparsity threshold, now with normalization $\sqrt{\E N_s}$.

\begin{theorem}
  \label{thm:upper-gaussian}
  If $d_s/\log s\to
  \beta\in(\beta_c,\beta_{\mathrm{sp}})$,
  then
  \begin{displaymath}
    \frac{N_s-\E N_s}{\sqrt{\E N_s}}
    \wto N(0,1).
  \end{displaymath}
  Moreover, if $\beta>\beta_1$, then $\Var N_s\sim\E N_s$, so
  $\sqrt{\E N_s}$ can be replaced by $\sqrt{\Var N_s}$. If
  $\beta<\beta_1$, then $\E N_s=o(\Var N_s)$ and thus
  \begin{displaymath}
    \frac{N_s-\E N_s}{\sqrt{\Var N_s}} \pto0.
  \end{displaymath}
\end{theorem}

We finally turn to the sparsity threshold
$\beta_{\mathrm{sp}}$. Assertions (ii) and (iii) below follow from
\cite{IlienkoJin2026} and are included for completeness.

\begin{theorem}
  \label{thm:densesp}~
  \begin{enumerate}[label=\textup{(\roman*)}]
  \item If $d_s/\log s\to\beta_{\mathrm{sp}}$
    and $d_s-\beta_{\mathrm{sp}}\log s\to-\infty$,
    then
    \begin{displaymath}
      \frac{N_s-\E N_s}{\sqrt{\Var N_s}}
      \wto N(0,1).
    \end{displaymath}
  \item \label{it:densesp-ii}
  If $d_{s_n}-\beta_{\mathrm{sp}}\log s_n
  \to x\in\R$ along some sequence $s_n\to\infty$,
  then
  \begin{displaymath}
  	N_{s_n}\wto\operatorname{Poisson}(2^{-x}).
  \end{displaymath}
  \item If $d_s-\beta_{\mathrm{sp}}\log s
    \to+\infty$, then $N_s\pto0$. In particular, this holds if
    $d_s/\log s\to\beta>\beta_{\mathrm{sp}}$.
  \end{enumerate}
\end{theorem}

The same remark as after
Theorem~\ref{thm:boundary-trichotomy} applies here.

The analogues of these results for a binomial process, that is a
sample of $n$ independent uniform points in $Q_{d_n}$, with
$d_n/\log n\to\beta>0$, are established in
Section~\ref{sec:de-poiss-binom} by a de-Poissonization argument; see
Theorem~\ref{thm:binomial-main} and
Corollary~\ref{cor:binomial-critical}.

\section{Martingale and Poisson chaos decompositions}
\label{sec:compensator-chaos}

In this section, we provide a decomposition of $N_s-\E N_s$ into a
compensated Poisson integral, the compensator first chaos, and a
higher-chaos remainder, and derive an exact variance identity.
Throughout this section, $s$ and $d=d_s$ are fixed. Write
$ \mu_s(\!\dd z)=\E\,\eta_s(\!\dd z)=s\dd z$ for the intensity measure
of $\eta_s$, and denote by $\widehat\eta_s=\eta_s-\mu_s$ the
compensated Poisson measure. For $k\geq1$, let
$L^2(\mu_s^{\otimes k})$ be the space of measurable kernels
$h:Q_d^k\to\R$ that are square-integrable with respect to the product
measure $\mu_s^{\otimes k}$. For symmetric
$h\in L^2(\mu_s^{\otimes k})$, let $I_k(h)$ denote the $k$-th multiple
Wiener--It\^o integral with respect to $\widehat\eta_s$, with the
normalization \begin{equation}\label{eq:orthchaos}
  \E[I_k(h)I_\ell(g)]
  =\one_{k=\ell} k!\,\langle h,g\rangle_{L^2(\mu_s^{\otimes k})},
\end{equation}
see \cite[Chapter~12]{LastPenrose}. In particular,
\begin{equation}\label{eq:iso}
  \E I_k(h)=0,
  \qquad
  \E I_k(h)^2
  =k!\,\|h\|_{L^2(\mu_s^{\otimes k})}^2,
  \qquad k\geq1.
\end{equation}
For $h\in L^2(\mu_s)$, the first compensated Poisson integral is
\begin{displaymath}
  I_1(h)=\int_{Q_d}h(z)\,\widehat\eta_s(\!\dd z)
  = \int_{Q_d}h(z)\,\eta_s(\!\dd z)
  -s\int_{Q_d}h(z)\dd z .
\end{displaymath}
In particular, from \eqref{eq:iso} we obtain
\begin{equation}\label{eq:VarI1}
  \Var I_1(h)
  =  \|h\|_{L^2(\mu_s)}^2
  =  s\int_{Q_d}h(z)^2\dd z.
\end{equation}

We next introduce the predictable decomposition of $N_s$ by
revealing points in increasing order of their coordinate product.
For $t\in[0,1]$, let
\begin{equation}
	\label{eq:Fst-def}
	\sF_{s,t} = \sigma\bigl(\eta_s(B):
	B\in\mathcal B(Q_d),\
	B\subseteq\{z\in Q_d:|z|\le t\}\bigr),
	\qquad 0\le t\le1.
\end{equation}
We complete the filtration $(\sF_{s,t})_{0\le t\le1}$,
which is then right-continuous. The number of non-Pareto points can be written as
\begin{displaymath}
  N_s=\int_{Q_d}H_{s,x}\,\eta_s(\!\dd x) = \sum_{x \in \eta_s} H_{s,x},
\end{displaymath}
where
\begin{equation}
  \label{eq:Hsx}
  H_{s,x}=\one_{\eta_s([0,x))\geq1} = \one\{x\text{ is not Pareto}\},
\end{equation}
and denote $\overline{H}_{s,x}=1-H_{s,x}$.
The jointly measurable map
$(x,\eta)\mapsto\one_{\{\eta([0,x))\ge1\}}$ satisfies
\begin{displaymath}
  H_{s,x}(\eta)
  =H_{s,x}\!\big(\eta|_{\{z:\,|z|<|x|\}}\big).
\end{displaymath}
This is exactly the predictability
condition for the ordered Poisson integral
\cite[Equation (2.5)]{LastPenroseMartingale}; it yields a predictable integrand for
$(\sF_{s,t})_{t\in[0,1]}$.

Define
\begin{equation}
  \label{eq:MsAs}
  M_s(t)=\int_{x\in Q_d,\,|x|\le t}
  H_{s,x}\,\widehat\eta_s(\!\dd x),
  \qquad
  A_s(t)=s\int_{x\in Q_d,\,|x|\le t}
  H_{s,x}\dd x,
\end{equation}
and set $M_s=M_s(1)$, $A_s=A_s(1)$.  By
\cite[Proposition~2.5]{LastPenroseMartingale}, $M_s(\cdot)$ is a
square-integrable martingale with predictable quadratic variation
$A_s(\cdot)$.  Moreover, $\E M_s=0$, $N_s=M_s+A_s$, and
\begin{equation}
  \label{eq:lambda-def}
  \E N_s  =
  \E A_s = s \int_{Q_d} \P\{\eta_s([0,x))\geq1\}\dd x
  =  s\int_{Q_d}\bigl(1-e^{-s|x|}\bigr)\dd x.
\end{equation}
Also, since $H_{s,x}^2=H_{s,x}$, by 
\cite[Equation (2.7)]{LastPenroseMartingale} we have
\begin{equation}
\label{eq:EM2}
  \E M_s^2
  =  \E\int_{Q_d}H_{s,x}^2\,\mu_s(\!\dd x)
  =  \E A_s
  =  \E N_s.
\end{equation}
Note that $1-A_s/s$ is the volume of the
undominated region $\{x\in Q_d:\eta_s([0,x))=0\}$.
Thus, limit theorems for $A_s$ immediately
translate into limit theorems for this volume.
As will be shown below, for $\beta\in(\beta_c,\beta_{\mathrm{sp}})$,
$A_s$ retains stable fluctuations, which are negligible on the scale
of the asymptotically Gaussian martingale $M_s$. Thus, in this regime
the volume of the undominated region exhibits stable fluctuations,
whereas $N_s$ has a Gaussian limit.

For a measurable Poisson functional $F=F(\eta_s)$,
define the Malliavin derivative and its iterates by
\[
  D_zF(\eta_s)=F(\eta_s+\delta_z)-F(\eta_s),\qquad
  D^k_{z_1,\ldots,z_k}F=D_{z_1}\cdots D_{z_k}F.
\]

\begin{proposition}
  \label{prop:chaosexp}  
  The centered compensator has the orthogonal expansion
  \begin{displaymath}
    A_s-\E A_s =
    \sum_{k=1}^\infty\frac1{k!}I_k(T_{k,s})
    \qquad\text{in }L^2(\P),
  \end{displaymath}
  where
  \begin{equation}
    \label{eq:ak-detailed}
    T_{k,s}(z_1,\ldots,z_k)
    = \E D^k_{z_1,\ldots,z_k}A_s =
    (-1)^{k+1}s
    \int_{z_1\vee\cdots\vee z_k}^1
    e^{-s|x|}\dd x,
  \end{equation}
  and $z_1\vee\cdots\vee z_k$ stands for the coordinatewise maximum.
\end{proposition}
\begin{proof}
  Since $0\le A_s\le s$, $A_s\in L^2(\P)$. By the Wiener--It\^o
  expansion \cite[Theorem~18.10]{LastPenrose},
  \begin{displaymath}
    A_s-\E A_s
    =
    \sum_{k=1}^\infty
    \frac1{k!}I_k(\E D^k_{\textstyle\cdot} A_s).
  \end{displaymath}
  Adding a point changes $H_{s,x}$ only if it lies in $[0,x)$
  and this set contains no point of $\eta_s$. Thus
  \begin{equation}
    \label{eq:DH}
    D_zH_{s,x}=\overline H_{s,x}\one_{\{z\in[0,x)\}},
    \qquad
    D_z\overline H_{s,x}=-\overline H_{s,x}\one_{\{z\in[0,x)\}}.
  \end{equation}
  Iterating this gives
  \begin{displaymath}
    D^k_{z_1,\ldots,z_k}H_{s,x}
    =
    (-1)^{k+1}\overline{H}_{s,x}
    \prod_{j=1}^k\one_{\{z_j\in[0,x)\}}.
  \end{displaymath}
  Since $\E\overline H_{s,x}=e^{-s|x|}$, integration yields
  \[
    \E D^k_{z_1,\ldots,z_k}A_s
    =(-1)^{k+1}s\int_{Q_d}e^{-s|x|}
    \prod_{j=1}^k\one_{\{z_j\in[0,x)\}}\dd x = T_{k,s}(z_1,\ldots,z_k).\qedhere
  \]
\end{proof}

By \eqref{eq:ak-detailed} and Fubini's theorem, for $k\ge1$ we have
\begin{equation}
  \label{eq:normT}
  \begin{aligned}
    \|T_{k,s}\|_{L^2(\mu_s^{\otimes k})}^2
    &=s^{k+2}\int_{Q_d^k}
    \biggl(\int_{z_1\vee\ldots\vee z_k}^1
    e^{-s|x|}\dd x\biggr)^2
    \dd z_1\ldots\dd z_k\\
    &=s^{k+2}\int_{Q_d^2}
    |x\wedge y|^k e^{-s(|x|+|y|)}
    \dd x\dd y,
  \end{aligned}
\end{equation}
where $x\wedge y$ denotes the coordinatewise minimum.
Taking $k=1$ in \eqref{eq:ak-detailed} and \eqref{eq:normT},
we obtain the compensator first-chaos kernel and its variance:
\begin{gather}
  \label{eq:fs-detailed}
  f_s(z)=T_{1,s}(z)=s\int_z^1 e^{-s|x|}\dd x,
  \qquad z\in Q_d,\\
  \label{eq:V1}
  V_{1,s}=\Var I_1(f_s)
  =s\int_{Q_d}f_s(z)^2\dd z
  =s^3\int_{Q_d^2}
  |x\wedge y|e^{-s(|x|+|y|)}\dd x\dd y.
\end{gather}
Proposition~\ref{prop:chaosexp} and \eqref{eq:orthchaos} yield
\begin{equation}
  \label{def:Rs-sum}
  \Var(A_s)=V_{1,s}+R_s,
  \qquad\text{where }
  R_s=\sum_{k=2}^\infty\frac1{k!}
  \|T_{k,s}\|_{L^2(\mu_s^{\otimes k})}^2.
\end{equation}
Summing the series using \eqref{eq:normT} gives
\begin{equation}
  \label{eq:R-def}
  R_s=s^2\int_{Q_d^2}e^{-s(|x|+|y|)}
  \left(e^{s|x\wedge y|}-1-s|x\wedge y|\right)
  \dd x\dd y.
\end{equation}

Throughout the paper we write $\sigma_s=\sqrt{V_{1,s}}$ and
\[
  X_s=\frac{I_1(f_s)}{\sigma_s}
  =\int_{Q_d}\frac{f_s(z)}{\sigma_s}
  \,\bigl(\eta_s(\!\dd z)-s\dd z\bigr)
\]
for the normalized compensator first chaos.  The Laplace functional of
a Poisson point process (see, e.g., \cite[Theorem~3.9]{LastPenrose})
gives
\begin{equation}
  \label{eq:Xs-chf}
  \E e^{itX_s}
  =\exp\left\{s\int_{Q_d}
    \left(e^{itf_s(z)/\sigma_s}-1
      -it\frac{f_s(z)}{\sigma_s}\right)\dd z
  \right\},\qquad t\in\R.
\end{equation}

By the Wiener--It\^o expansion, the first chaos
of $N_s-\E N_s$ is $I_1(g_s)$, where
\begin{equation}
  \label{eq:gs}
  g_s(z)=\E D_zN_s=f_s(z)+\bigl(1-e^{-s|z|}\bigr)=f_s(z)+q_s(z).
\end{equation}
Indeed, $D_zN_s=N_s(\eta_s+\delta_z)-N_s(\eta_s)$
counts the added point if it is non-Pareto, together with
the Pareto points that become non-Pareto after its insertion.
The first contribution has expectation $q_s(z)$ while the
second has expectation
$s\int_z^1e^{-s|x|}\dd x=f_s(z)$.
As $M_s=N_s-A_s$,
the first chaoses of $A_s-\E A_s$ and $M_s$
are $I_1(f_s)$ and $I_1(q_s)$, respectively.
Throughout, we call $I_1(f_s)$ the
compensator first chaos and $I_1(g_s)$ the full
first chaos. Set
\begin{equation}
  \label{eq:Us}
  U_s=A_s-\E A_s-I_1(f_s).
\end{equation}
Thus $U_s$ is the higher-chaos remainder of the compensator.
By \eqref{def:Rs-sum} and orthogonality, we have
$\E U_s^2=R_s$, and
\begin{equation}
  \label{eq:N-EN-expression}
  N_s-\E N_s=M_s+(A_s-\E A_s)=I_1(f_s)+M_s+U_s.
\end{equation}
Consequently, by the triangle inequality and \eqref{eq:EM2},
\begin{equation}
  \label{eq:bound-difference-chaos}
  \|N_s-\E N_s-I_1(f_s)\|_2
  \le\sqrt{\E M_s^2}+\sqrt{\E U_s^2}
  =\sqrt{\E N_s}+\sqrt{R_s}.
\end{equation}

Finally, we establish the orthogonality of $M_s$ and $A_s$
and derive an exact variance identity.

\begin{proposition}
  \label{prop:variance-identity}
  The following identities hold:
  \[
    \Cov(M_s,A_s)=0,
    \qquad
    \Var N_s=\E N_s+\Var A_s
    =\E N_s+V_{1,s}+R_s.
  \]
\end{proposition}
\begin{proof}
  Since $M_s=\int_{Q_d}H_{s,z}\,\widehat\eta_s(\!\dd z)$, we have
  \[
    \E[M_sA_s]
    =\E\int_{Q_d}H_{s,z}A_s\,\eta_s(\!\dd z)
    -s\int_{Q_d}\E[H_{s,z}A_s]\dd z.
  \]
  By the Mecke formula and the fact that
  $H_{s,z}(\eta_s+\delta_z)=H_{s,z}(\eta_s)$,
  \[
    \E\int_{Q_d}H_{s,z}A_s\,\eta_s(\!\dd z)
    =s\int_{Q_d}\E[H_{s,z}A_s^{+z}]\dd z,
  \]
  where $A_s^{+z}$ denotes the compensator after adding
  a point at $z$. Since $A_s^{+z}-A_s=D_zA_s$,
  \[
    \E[M_sA_s]
    =s\int_{Q_d}\E[H_{s,z}D_zA_s]\dd z
    =s^2\int_{Q_d}\int_z^1
    \E[H_{s,z}\overline H_{s,x}]
    \dd x\dd z=0,
  \]
  where the second equality uses
  $D_zA_s=s\int_z^1\overline H_{s,x}\dd x$ from \eqref{eq:DH}, and the
  last uses $H_{s,z}\overline H_{s,x}=0$ whenever $z\le x$.  As
  $\E M_s=0$, this proves $\Cov(M_s,A_s)=0$.  The variance identity
  follows from $N_s=M_s+A_s$, \eqref{eq:EM2}, and \eqref{def:Rs-sum}.
\end{proof}

\section{Asymptotics for moment integrals}
\label{sec:prod-minim-lapl}

In this section, we determine the asymptotics of a class of
integrals used in our variance and moment estimates. For fixed
$m\in\mathbb{N}$, $\rho\ge0$, and $c_1,\ldots,c_m>0$, define
\begin{equation}
  \label{eq:J-general}
  J^{(c)}_{m,\rho}(s,d_s)
  =\int_{(Q_d)^m}|z_1\wedge\cdots\wedge z_m|^\rho
  \exp\biggl\{-s\sum_{r=1}^m c_r|z_r|\biggr\}
  \dd z_1\cdots\dd z_m,\qquad s>0,\quad d_s\in\mathbb{N}.
\end{equation}
We shall also use
\begin{equation}
  \label{eq:M}
  M_{m,\rho}(t)=\frac{m!}{\prod_{j=1}^m(\rho+j-jt)},
  \qquad 0\le t<1+\frac\rho m,
\end{equation}
whose meaning will become clear in the proof. The following
proposition gives the required asymptotics.

\begin{proposition}
  \label{prop:laplace}
  Let $d_s/\log s\to\beta\in(0,\infty)$. Then
  \begin{equation}
    \label{eq:general-exponent}
    J^{(c)}_{m,\rho}(s,d_s)=s^{\xi_{m,\rho}(\beta)+o(1)},
  \end{equation}
  where
  \begin{equation}
    \label{eq:xi-variational}
    \xi_{m,\rho}(\beta)
    =\min_{0\le t<1+\rho/m}
    \bigl\{\beta\log M_{m,\rho}(t)-mt\bigr\}.
  \end{equation}
  For fixed $m$ and $\rho$, the minimizer $t_\beta$ is unique and satisfies
  \[
    t_\beta>0\quad\text{if }\beta<\beta_{m,\rho}^{\rm c}
    \qquad\text{and}\qquad
    t_\beta=0\quad\text{if }\beta\ge\beta_{m,\rho}^{\rm c},
  \]
  where
  \begin{equation}
    \label{eq:betamrhoc}
    \beta_{m,\rho}^{\rm c}
    =\biggl(1-\frac{\rho}{m}\sum_{j=1}^m\frac1{\rho+j}\biggr)^{-1}.
  \end{equation}
\end{proposition}
\begin{proof}
  Let $B_{ri}$, $1\le r\le m$, $1\le i\le d_s$, be independent
  $\operatorname{Exp}(1)$ random variables, and put
  \[
    \Gamma_r=\sum_{i=1}^{d_s}B_{ri},\qquad
    R=\sum_{i=1}^{d_s}(B_{1i}\vee\cdots\vee B_{mi}).
  \]
  The logarithmic change of variables $z_{ri}=e^{-b_{ri}}$ yields
  \begin{equation}
    \label{eq:Jc}
    J^{(c)}_{m,\rho}(s,d_s)
    =\E\exp\biggl\{-\rho R-s\sum_{r=1}^m c_re^{-\Gamma_r}\biggr\}.
  \end{equation}
  To analyze \eqref{eq:Jc}, we first consider a simpler exponential
  moment:
  \[
    \E\exp\biggl\{t\sum_{r=1}^m\Gamma_r-\rho R\biggr\}
    =\biggl(\E\exp\biggl\{t\sum_{r=1}^mB_{r1}
    -\rho(B_{11}\vee\cdots\vee B_{m1})\biggr\}\biggr)^{d_s}
    =M_{m,\rho}(t)^{d_s},
  \]
  where the second equality readily follows from the standard
  independent-spacing representation for exponential order statistics.
  By \eqref{eq:M}, $\beta\log M_{m,\rho}(t)-mt$ is strictly convex in
  $t$ and tends to infinity as $t\uparrow1+\rho/m$. Hence it has a
  unique minimizer $t_\beta$.
	
  Let $\P$ be the law of $(B_{11},\ldots,B_{m1})$ and define $\P_t$ by
  \begin{equation}
    \label{eq:PtP}
    \frac{\dd\P_t}{\dd\P}
    =\frac1{M_{m,\rho}(t)}
    \exp\biggl\{t\sum_{r=1}^mB_{r1}
    -\rho(B_{11}\vee\cdots\vee B_{m1})\biggr\}.
  \end{equation}
  Applying this change of measure independently to the $d_s$
  coordinates and writing $\E_t$ for expectation under the resulting
  law, we obtain from \eqref{eq:Jc}
  \begin{equation}
    \label{eq:J-under-tilt-section4}
    \begin{aligned}
      J^{(c)}_{m,\rho}(s,d_s)
      &=M_{m,\rho}(t)^{d_s}\,\E_t\exp\biggl\{-t\sum_{r=1}^m\Gamma_r
      -s\sum_{r=1}^m c_re^{-\Gamma_r}\biggr\}\\
      &=M_{m,\rho}(t)^{d_s} s^{-mt}\,
      \E_t\prod_{r=1}^m\exp\{-tX_r-c_re^{-X_r}\},
    \end{aligned}
  \end{equation}
  where $X_r=\Gamma_r-\log s$.
  
  Since, for $t\ge0$ and $c>0$, $\exp\{-tx-ce^{-x}\}$ is bounded on
  $\R$, \eqref{eq:J-under-tilt-section4} gives
  \[
    J^{(c)}_{m,\rho}(s,d_s)\le C M_{m,\rho}(t)^{d_s} s^{-mt},
  \]
  where $C$ does not depend on $s$. Therefore, for every
  $t\in[0,1+\rho/m)$,
  \[
    \limsup_{s\to\infty}
    \frac{\log J^{(c)}_{m,\rho}(s,d_s)}{\log s}
    \le\beta\log M_{m,\rho}(t)-mt.
  \]
  Thus, the left-hand side is bounded above by $\xi_{m,\rho}(\beta)$.
  
  For the reverse bound, put
  $\mu(t)=\frac1m\frac{\dd}{\dd t}\log M_{m,\rho}(t)$. Note that, since
  $\frac{\dd\P_t}{\dd\P}$ in \eqref{eq:PtP} is symmetric in
  $B_{11},\ldots,B_{m1}$, we also have
  \begin{align*}
    \E_tB_{11}
    &=\frac1m\E_t\sum_{r=1}^mB_{r1}\\
    &=\frac1m
      \frac{\E\bigl[(\sum_{r=1}^mB_{r1})
      \exp\{t\sum_{r=1}^mB_{r1}-\rho(B_{11}\vee\cdots\vee B_{m1})\}\bigr]}
      {\E\exp\{t\sum_{r=1}^mB_{r1}-\rho(B_{11}\vee\cdots\vee B_{m1})\}}
      =\frac1m\frac{\dd M_{m,\rho}(t)/\!\dd t}{M_{m,\rho}(t)}
      =\mu(t).
  \end{align*}
  
  Note that in both cases $t_\beta>0$ and $t_\beta=0$ we have
  $\mu(t_\beta)\ge\beta^{-1}$, with equality if $t_\beta>0$. Since
  $\log M_{m,\rho}$ is strictly convex, $\mu(t)>\beta^{-1}$ for all
  $t\in(t_\beta,1+\rho/m)$. Fix such a $t$ and choose
  $\eps\in(0,\mu(t)-\beta^{-1})$. By the law of large numbers,
  \[
    E_s=\bigcap_{r=1}^m\{d_s(\mu(t)-\eps)\le \Gamma_r\le d_s(\mu(t)+\eps)\}
  \]
  satisfies $\P_t(E_s)\to1$. On $E_s$, we have
  \[
    \sum_{r=1}^m\Gamma_r\le md_s(\mu(t)+\eps)\qquad\text{and}\qquad
    s\sum_{r=1}^m c_re^{-\Gamma_r}
    \le\Bigl(\sum_{r=1}^m c_r\Bigr)se^{-d_s(\mu(t)-\eps)}\to0,
  \]
  since $\log s/d_s\to\beta^{-1}<\mu(t)-\eps$. Restricting the first expectation in
  \eqref{eq:J-under-tilt-section4} to $E_s$, we obtain
  \[
    J^{(c)}_{m,\rho}(s,d_s)
    \ge M_{m,\rho}(t)^{d_s}\P_t(E_s)
    \exp\{-mtd_s(\mu(t)+\eps)-o(1)\}.
  \]
  Taking logarithms, dividing by $\log s$, and letting $\eps\downarrow0$
  gives
  \begin{equation*}
    \liminf_{s\to\infty}
    \frac{\log J^{(c)}_{m,\rho}(s,d_s)}{\log s}
    \ge\beta\{\log M_{m,\rho}(t)-mt\mu(t)\}.
  \end{equation*}
  Letting $t\downarrow t_\beta$, we obtain
  \[
    \liminf_{s\to\infty}\frac{\log J^{(c)}_{m,\rho}(s,d_s)}{\log s}
    \ge\beta\log M_{m,\rho}(t_\beta)
    -\beta mt_\beta\mu(t_\beta)=\xi_{m,\rho}(\beta),
  \]
  where the last equality follows from
  $\mu(t_\beta)=\beta^{-1}$ if $t_\beta>0$ and is immediate if
  $t_\beta=0$.	
  Together with the upper bound, this proves \eqref{eq:general-exponent}.
  
  Finally, it readily follows from \eqref{eq:M} that
  \[
    \mu(0)=1-\frac{\rho}{m}\sum_{j=1}^m\frac1{\rho+j}.
  \]
  Since the right derivative at $0$ of $\beta\log M_{m,\rho}(t)-mt$
  has the same sign as $\beta\mu(0)-1$, strict convexity shows that
  $t_\beta>0$ precisely when $\beta<1/\mu(0)=\beta_{m,\rho}^{\rm
    c}$. This completes the proof.
\end{proof}

For $\beta>\beta_{m,\rho}^{\rm c}$, the asymptotics can be sharpened
and extended to more general integrands.

\begin{proposition}
  \label{cor:strong-zero-saddle}
  Let $d_s/\log s\to\beta>\beta_{m,\rho}^{\rm c}$. If
  $\psi_1,\ldots,\psi_m:[0,\infty)\to\R$ are bounded
  measurable functions such that each $\psi_r$ is continuous at $0$
  and $\psi_r(0)=1$, then
  \begin{equation}
    \label{eq:strong-zero-saddle-weighted}
    \int_{(Q_{d_s})^m}|z_1\wedge\cdots\wedge z_m|^\rho
    \prod_{r=1}^m\psi_r\bigl(sc_r|z_r|\bigr)
    \,\dd z_1\cdots\dd z_m
    \sim M_{m,\rho}(0)^{d_s}.
  \end{equation}
\end{proposition}
\begin{proof}
  For $\beta>\beta_{m,\rho}^{\rm c}$, we have $t_\beta=0$ and
  $\mu(0)>\beta^{-1}$. Hence, similarly to the first equality in
  \eqref{eq:J-under-tilt-section4}, the left-hand side of
  \eqref{eq:strong-zero-saddle-weighted} equals
  \[
    M_{m,\rho}(0)^{d_s}
    \E_0\prod_{r=1}^m\psi_r\bigl(sc_re^{-\Gamma_r}\bigr).
  \]
  By the law of large numbers,
  \[
    \frac{\Gamma_r}{\log s}
    =\frac{\Gamma_r}{d_s}\frac{d_s}{\log s}
    \xrightarrow{\P_0}\beta\,\E_0B_{11}
    =\beta\mu(0)>1.
  \]
  Hence $sc_re^{-\Gamma_r}\to0$ in $\P_0$-probability for each $r$.
  The product in the expectation is bounded and converges to $1$ in
  probability, so its expectation tends to $1$.
\end{proof}

\section{Asymptotics of the variance components}
\label{sec:comp-first-chaos}

Proposition~\ref{prop:variance-identity} shows that the variance of the number $N_s$ of non-Pareto points decomposes into its mean $\E N_s$, the variance $V_{1,s}$ of the first-chaos term of the centered compensator, and the residual variance $R_s$. In this section, we compare these three terms. In particular, we identify the range where the compensator first chaos dominates and derive several asymptotics used below.

For $m=2$, $\rho=1$, \eqref{eq:M} and \eqref{eq:betamrhoc} give
\[
M_{2,1}(t)=\frac{2}{(2-t)(3-2t)}, \qquad 0\le t<\frac32, \qquad\text{and}\qquad\beta_{2,1}^{\rm c}=\frac{12}{7}.
\]
For $\beta\in(0,12/7)$, the minimizer $t_\beta$ in \eqref{eq:xi-variational} is positive and satisfies
\[
\frac12\left(\frac1{2-t_\beta}+\frac2{3-2t_\beta}\right)=\frac1\beta.
\]
Thus
\begin{equation}\label{eq:tb}
	t_\beta=\frac{7-2\beta-\sqrt{1+4\beta^2}}4.
\end{equation}
For $\beta\ge12/7$, we have $t_\beta=0$.
Hence \eqref{eq:V1}, \eqref{eq:J-general}, and Proposition~\ref{prop:laplace} yield
\begin{equation}\label{eq:V1-power}
	V_{1,s}=s^3J^{(1,1)}_{2,1}(s,d_s)=s^{r(\beta)+o(1)},
\end{equation}
where $r(\beta)=3+\xi_{2,1}(\beta)$. Explicitly,
\begin{equation}
	\label{eq:rbeta-saddle}
	r(\beta)=
	\begin{cases}
		3-2t_\beta+
		\beta\log\frac{2}{(2-t_\beta)(3-2t_\beta)},
		&0<\beta<12/7,\\
		3-\beta\log3,
		&\beta\ge12/7.
	\end{cases}
\end{equation}
Let $\beta_0\approx0.305$ be the smaller positive solution of $r(\beta)=1$.
For $\beta>12/7$, Proposition~\ref{cor:strong-zero-saddle} yields the stronger asymptotic
\begin{equation}
	\label{eq:V1s}
	V_{1,s}\sim s^3M_{2,1}(0)^{d_s}=s^33^{-d_s}.
\end{equation}

We next estimate the higher-order chaos remainder.

\begin{lemma}
	\label{lem:R-bounds}
	Let $d_s/\log s\to\beta>0$. Then $R_s\le\frac{s^4}26^{-d_s}$.
	Moreover, if $\beta\le12/7$, then $R_s=o(V_{1,s})$.
\end{lemma}
\begin{proof}
	From \eqref{eq:R-def}, using $e^u-1-u\le u^2e^u/2$ and $|x\wedge y|\le(|x|+|y|)/2$, we obtain
	\begin{equation}\label{eq:Rs-bound-1}
		R_s\le\frac{s^4}{2}\int_{Q_{d_s}^2}|x\wedge y|^2e^{-s(|x|+|y|)/2}\dd x\dd y\le
		\frac{s^4}{2}\int_{Q_{d_s}^2}|x\wedge y|^2\dd x\dd y=\frac{s^4}{2}6^{-d_s},
	\end{equation}
	which proves the first claim. For the second claim, let $t_\beta$ be the minimizer in \eqref{eq:xi-variational} for $(m,\rho)=(2,1)$, explicitly given by \eqref{eq:tb}. Since $0\le t_\beta<3/2$ for $\beta\in(0,12/7]$, we have $t_\beta+1/2<2$, so this is an admissible value in \eqref{eq:xi-variational} for $(m,\rho)=(2,2)$. Thus, by \eqref{eq:Rs-bound-1} and Proposition~\ref{prop:laplace},
	\begin{equation*}
		R_s\le\frac{s^4}{2}J_{2,2}^{(1/2,1/2)}(s,d_s)
		=s^{4+\xi_{2,2}(\beta)+o(1)}
		\le s^{4+\beta\log M_{2,2}(t_\beta+1/2)-2(t_\beta+1/2)+o(1)}.
	\end{equation*}
	Moreover,
	\[
	M_{2,2}(t_\beta+1/2)=\frac{2}{(3-(t_\beta+1/2))(4-2(t_\beta+1/2))}
	<\frac{2}{(2-t_\beta)(3-2t_\beta)}=M_{2,1}(t_\beta).
	\]
	Therefore,
	\[
	4+\beta\log M_{2,2}(t_\beta+1/2)-2(t_\beta+1/2)
	<3+\beta\log M_{2,1}(t_\beta)-2t_\beta=r(\beta).
	\]
	Together with \eqref{eq:V1-power}, this yields $R_s=o(V_{1,s})$.
\end{proof}

\begin{proposition}
	\label{prop:variance-dominance}
	For $\beta\in(\beta_0,\beta_1)$, $\Var N_s\sim V_{1,s}$.
\end{proposition}

\begin{proof}
	By Proposition~\ref{prop:variance-identity},
	\[
	\Var N_s=\E N_s+V_{1,s}+R_s.
	\]
	We first show that $\E N_s=o(V_{1,s})$. For $\beta\in(\beta_0,12/7]$, by \eqref{eq:xi-variational},
	\[
	r(\beta)=3+\inf_{0\le t<3/2}\{\beta\log M_{2,1}(t)-2t\}.
	\]
	Thus $r$ is concave as the infimum of affine functions of $\beta$. Since $r(\beta_0)=1$ by definition of $\beta_0$ and $r(12/7)=3-\frac{12}{7}\log3>1$,
	we have $r(\beta)>1$ for $\beta\in(\beta_0,12/7]$. Hence, using $\E N_s\le s$ and \eqref{eq:V1-power}, we obtain $\E N_s=o(V_{1,s})$.
	If instead $\beta\in(12/7,\beta_1)$, then
	from \eqref{eq:lambda-def} and $1-e^{-u}\le u$,
	\begin{equation*}
		\E N_s\le s^2\int_{Q_{d_s}}|x|\dd x=s^22^{-d_s}.
	\end{equation*}
	Since $V_{1,s}\sim s^33^{-d_s}$ for $\beta>12/7$,
	\begin{equation}
	\label{eq:NV}
	\frac{\E N_s}{V_{1,s}}\le s^{-1}\left(\frac32\right)^{d_s}(1+o(1))
	=s^{-1+\beta\log(3/2)+o(1)}\to0,
	\end{equation}
	because $\beta<\beta_1=1/\log(3/2)$.
	
	It remains to show that $R_s=o(V_{1,s})$. This follows from Lemma~\ref{lem:R-bounds} for $\beta\le12/7$. For $\beta\in(12/7,\beta_1)$, the same lemma and $V_{1,s}\sim s^33^{-d_s}$ give
	\[
	\frac{R_s}{V_{1,s}}\le\frac12s\,2^{-d_s}(1+o(1))
	=s^{1-\beta\log2+o(1)}\to0,
	\]
	since $\beta>12/7>1/\log2$.
\end{proof}

Below, for a sequence of random variables $(Y_s)$ and a positive sequence $(a_s)$, we write $Y_s=o_{L^2}(a_s)$ if $Y_s/a_s\to0$ in $L^2$. Recall that $\sigma_s=\sqrt{V_{1,s}}$, $f_s(z)=s\int_{[z,1]}e^{-s|x|}\dd x$, and $g_s(z)=f_s(z)+1-e^{-s|z|}$.

\begin{corollary}
	\label{cor:first-chaos}
	For $\beta\in(\beta_0,\beta_1)$,
	\[
	M_s=o_{L^2}(\sigma_s),\qquad
	A_s-\E A_s=I_1(f_s)+o_{L^2}(\sigma_s),\qquad
	N_s-\E N_s=I_1(f_s)+o_{L^2}(\sigma_s).
	\]
	Moreover, the full first-chaos projection $I_1(g_s)$ of $N_s-\E N_s$ satisfies
	$I_1(g_s)-I_1(f_s)=o_{L^2}(\sigma_s)$.
\end{corollary}

\begin{proof}
	The proof of Proposition~\ref{prop:variance-dominance} gives $\E N_s=o(V_{1,s})$ and $R_s=o(V_{1,s})$. Hence, by \eqref{eq:EM2}, $\|M_s\|_2^2=\E N_s=o(V_{1,s})$,
	while the orthogonality of the compensator chaoses and \eqref{def:Rs-sum} yield
	\[
	\|A_s-\E A_s-I_1(f_s)\|_2^2=R_s=o(V_{1,s}).
	\]
	Since $N_s-\E N_s=M_s+(A_s-\E A_s)$, the first three assertions follow.
	
	Finally, by \eqref{eq:gs}, $g_s=f_s+q_s$, where $q_s(z)=1-e^{-s|z|}$. Thus
	\[
	\|I_1(g_s)-I_1(f_s)\|_2^2=s\int_{Q_{d_s}}q_s(z)^2\,\dd z\le s\int_{Q_{d_s}}q_s(z)\,\dd z=\E N_s=o(V_{1,s}),
	\]
	which proves the last assertion.
\end{proof}

The following corollary shows that for $\beta\in(\beta_0,\beta_1)$,
any distributional limit of the normalized compensator first chaos
is also a limit of the centered count normalized by its standard
deviation. This follows immediately from the third assertion of
Corollary~\ref{cor:first-chaos} and
Proposition~\ref{prop:variance-dominance}.

\begin{corollary}
	\label{cor:transfer}
	Let $\beta\in(\beta_0,\beta_1)$. If $I_1(f_s)/\sigma_s\wto Y$ for some random variable $Y$\!, then
	\[
	\frac{N_s-\E N_s}{\sqrt{\Var N_s}}\wto Y.
	\]
\end{corollary}

We finish with a simple asymptotic estimate for $\E N_s$ that will be used below.

\begin{lemma}
	\label{lem:mean-strong}
	If $\beta>2$, then $\E N_s\sim s^22^{-d_s}$.
\end{lemma}

\begin{proof}
	By \eqref{eq:lambda-def},
	\[
	\E N_s=s^2\int_{Q_{d_s}}|x|\psi(s|x|)\,\dd x,
	\]
	where $\psi(u)=(1-e^{-u})/u$ for $u>0$
	and $\psi(0)=1$.
	Since $M_{1,1}(0)=1/2$ and $\beta_{1,1}^{\rm c}=2$, Proposition \ref{cor:strong-zero-saddle} yields the claim.
\end{proof}

\section{The Gaussian and boundary regime $0\leq \beta\leq\beta_L$}
\label{sec:gauss-rang-beta_0b}

In this section we prove Theorem~\ref{thm:main-gaussian} and the
boundary law at $\beta_L$. By Corollary~\ref{cor:transfer}, for
$\beta\in(\beta_0,\beta_L]$ it suffices to analyze the normalized first
chaos $X_s=I_1(f_s)/\sigma_s$. We establish convergence of $X_s$ to a centered Gaussian law, with variance one for $\beta\in(\beta_0,\beta_L)$ and possibly smaller variance at $\beta=\beta_L$, depending on the second-order behavior of the dimension.  The required estimates are established in
Sections~\ref{sec:positive-saddle} and~\ref{sec:zero-saddle}. We begin,
however, with the range $\beta\in[0,\beta_0]$, where a standard
Gaussian approximation still holds, but the fluctuation of the total
count $\eta_s(Q_{d_s})$ can no longer be neglected.

\subsection{The Gaussian range
	$0\leq\beta\leq\beta_0$}
\label{sec:appr-with-total}

For $\beta\le\beta_0$ we have $r(\beta)\le1$, so by \eqref{eq:V1-power}
$V_{1,s}=s^{r(\beta)+o(1)}$ need not dominate
$\Var M_s=\E N_s\le s$. We therefore approximate $N_s-\E N_s$ by its
full first-chaos projection $I_1(g_s)$, where $g_s=f_s+q_s$ and
$q_s(z)=1-e^{-s|z|}$ by \eqref{eq:gs}.
We will need the following Lyapunov-type lemma.

\begin{lemma}
	\label{lem:poisson-lyapunov}
	Let $h_s$ be a deterministic non-negative function such that
	$\int_{Q_{d_s}}h_s(z)^3\dd z<\infty$, and put
	$\zeta_s^2=s\int_{Q_{d_s}} h_s(z)^2\,\dd z>0$. If
	\[
	\frac s{\zeta_s^3}\int_{Q_{d_s}} h_s(z)^3\dd z\to0,
	\]
	then
	$I_1(h_s)/\zeta_s\wto N(0,1)$.
\end{lemma}
\begin{proof}
	Since $|e^{iu}-1-iu+u^2/2|\le|u|^3/6$, as in \eqref{eq:Xs-chf}, we have
	\begin{displaymath}
		\biggl|\log\E e^{itI_1(h_s)/\zeta_s}+\frac{t^2}{2}\biggr|
		=\biggl|s\int_{Q_{d_s}}\biggl(e^{ith_s(z)/\zeta_s}-1-\frac{ith_s(z)}{\zeta_s}
		+\frac{t^2h_s(z)^2}{2\zeta_s^2}\biggr)\dd z\biggr|
		\le\frac{|t|^3}{6} \frac s{\zeta_s^3}\int_{Q_{d_s}}h_s(z)^3\dd z,
	\end{displaymath}
	which tends to $0$ by assumption.
\end{proof}

\begin{proposition}
	\label{prop:low-gaussian}
	Suppose that $d_s/\log s\to\beta\in[0,\beta_0]$. Then
	\[
	\frac{N_s-\E N_s}{\sqrt{\Var N_s}}\wto N(0,1).
	\]
\end{proposition}

\begin{proof}
	Using the bound	$e^{-u}\le C_pu^{-p}$, $u>0$,
	with $p=2/3$, $9/5$, and $8/7$,
	we easily obtain
	\begin{gather}
		\notag
		s\int_{Q_{d_s}}e^{-s|z|}\dd z=o(s),
		\qquad
		R_s\le\frac{s^4}{2}\int_{Q_{d_s}^2}|x\wedge y|^2
		e^{-s(|x|+|y|)/2}\dd x\dd y=o(s),\\
		s\int_{Q_{d_s}}f_s(z)^3\dd z
		=s^4\iiint_{Q_{d_s}^3}|x\wedge y\wedge z|
		e^{-s(|x|+|y|+|z|)}\dd x\dd y\dd z
		=o(s^{3/2}).
		\label{eq:low-f-third}
	\end{gather}
	
	Since $I_1(q_s)$ is the first-chaos projection of $M_s$,
	\eqref{eq:EM2}, \eqref{eq:lambda-def}, and \eqref{eq:iso} yield
	\[
	\|M_s-I_1(q_s)\|_2^2=\|M_s\|_2^2-\|I_1(q_s)\|_2^2
	=s\int_{Q_{d_s}}q_s(z)(1-q_s(z))\dd z
	\le s\int_{Q_{d_s}}e^{-s|z|}\dd z=o(s).
	\]
	Moreover,
	\[
	\|A_s-\E A_s-I_1(f_s)\|_2^2=R_s=o(s).
	\]
	Hence, by \eqref{eq:N-EN-expression},
	\begin{equation}
		\label{eq:low-first-chaos}
		\|N_s-\E N_s-I_1(g_s)\|_2=o(\sqrt{s}).
	\end{equation}
	
	Put
	\[
	\tau_s^2=s\int_{Q_{d_s}}g_s(z)^2\dd z.
	\]
	Since $g_s\ge q_s$ and $1-q_s=e^{-s|z|}$,
	\[
	\tau_s^2
	\ge s\int_{Q_{d_s}}q_s(z)^2\dd z
	\ge s-2s\int_{Q_{d_s}}e^{-s|z|}\dd z
	=s-o(s).
	\]
	Hence, since $(a+b)^3\le4(a^3+b^3)$ and $0\le q_s\le1$,
	\[
	s\int_{Q_{d_s}}g_s(z)^3\dd z
	\le4s\int_{Q_{d_s}}f_s(z)^3\dd z+4s
	=o(s^{3/2})
	=o(\tau_s^3).
	\]
	Thus, Lemma~\ref{lem:poisson-lyapunov} yields
	\[
	\frac{I_1(g_s)}{\tau_s}\wto N(0,1).
	\]
	
	Since $I_1(g_s)$ is the first-chaos projection of $N_s-\E N_s$,
	\eqref{eq:low-first-chaos} and orthogonality give
	\[
	\Var N_s
	=\tau_s^2+\|N_s-\E N_s-I_1(g_s)\|_2^2
	=\tau_s^2+o(s)
	\sim\tau_s^2.
	\]
	Therefore, again by \eqref{eq:low-first-chaos},
	\[
	\frac{N_s-\E N_s}{\sqrt{\Var N_s}}
	=\frac{I_1(g_s)}{\tau_s}
	\frac{\tau_s}{\sqrt{\Var N_s}}+o_{\P}(1)
	\wto N(0,1).\qedhere
	\]
\end{proof}

\subsection{The Gaussian range $\beta_0<\beta\le12/7$}
\label{sec:positive-saddle}

By Corollary~\ref{cor:transfer}, it suffices to prove asymptotic
normality of $I_1(f_s)$. We again verify the Lyapunov condition from
Lemma~\ref{lem:poisson-lyapunov}. By the first equality in
\eqref{eq:low-f-third},
\[
s\int_{Q_{d_s}}f_s(z)^3\dd z=s^4J^{(1,1,1)}_{3,1}(s,d_s).
\]
Hence Proposition~\ref{prop:laplace} gives
\begin{equation}
	\label{eq:r3-def}
	s\int_{Q_{d_s}}f_s(z)^3\dd z=s^{r_3(\beta)+o(1)},\qquad
	\text{where }	r_3(\beta)=
	4+\min_{0\le v<4/3}
	\left\{\beta\log M_{3,1}(v)-3v\right\}.
\end{equation}
Since $\sigma_s^2=V_{1,s}=s^{r(\beta)+o(1)}$ by
\eqref{eq:V1-power}, it remains to show that
$r_3(\beta)<3r(\beta)/2$.

\begin{lemma}
	\label{lem:r3}
	For every $\beta\in(0,12/7]$, $r_3(\beta)<\frac32r(\beta)$.
\end{lemma}

\begin{proof}
	Let first $\beta\in(0,88/57]$. The saddle $t=t_\beta$ in \eqref{eq:tb}
	decreases with $\beta$ and satisfies $t_{88/57}=1/6$. Thus
	$v=t-1/6\in[0,4/3)$ is admissible in \eqref{eq:r3-def}, and
	\[
	r_3(\beta)\le
	4+\beta\log M_{3,1}(t-1/6)-3t+\frac12.
	\]
	Since
	\[
	r(\beta)=3+\beta\log M_{2,1}(t)-2t,
	\]
	we obtain
	\[
	r_3(\beta)-\frac32r(\beta)
	\le\frac{\beta}{2}
	\log\frac{M_{3,1}(t-1/6)^2}{M_{2,1}(t)^3}.
	\]
	Put $x=12-6t$. Then $3<x\le11$, and \eqref{eq:M} readily gives
	\[
	\frac{M_{3,1}(t-1/6)^2}{M_{2,1}(t)^3}
	=\frac{x^3(x-3)}{(x+1)^2(x-2)^2}<1.
	\]
	This proves the claim for $\beta\in(0,88/57]$.
	
	Let now $\beta\in(88/57,12/7]$. Taking $v=0$ in
	\eqref{eq:r3-def} gives $r_3(\beta)\le4-\beta\log4$.
	It therefore suffices to show that
	\[
	H(\beta):=4-\beta\log4-\frac32r(\beta)<0.
	\]
	Since $r$ is concave, $H$ is convex. Using
	$t_{88/57}=1/6$, $t_{12/7}=0$, and
	\eqref{eq:rbeta-saddle}, we obtain
	\[
	H(88/57)
	=\frac{88}{57}
	\log\left(\frac14\left(\frac{22}{9}\right)^{3/2}\right)<0,
	\qquad
	H(12/7)
	=-\frac12+\frac67\log\frac{27}{16}<0.
	\]
	Hence $H(\beta)<0$ throughout the interval.
\end{proof}

\begin{proposition}
	\label{prop:positive-gaussian}
	If $\beta\in(\beta_0,12/7]$, then
	\[
	\frac{N_s-\E N_s}{\sqrt{\Var N_s}}\wto N(0,1).
	\]
\end{proposition}

\begin{proof}
	By Lemma~\ref{lem:r3}, the Lyapunov condition in
	Lemma~\ref{lem:poisson-lyapunov} holds. Hence
	\[
	X_s=\frac{I_1(f_s)}{\sigma_s}\wto N(0,1),
	\]
	and the claim follows from Corollary~\ref{cor:transfer}.
\end{proof}

\subsection{The Gaussian range $12/7<\beta\le\beta_L$}
\label{sec:zero-saddle}

For $\beta>12/7$, \eqref{eq:V1s} gives the exact
asymptotic $V_{1,s}\sim s^33^{-d_s}$, which will allow us to replace $f_s$ by a simpler product kernel. Denote
\begin{equation}
	\label{eq:v-T-def}
	v(z)=\prod_{i=1}^{d_s}(1-z_i),
	\qquad
	T(z)=-\log v(z)=\sum_{i=1}^{d_s}\bigl(-\log(1-z_i)\bigr),
	\qquad z\in Q_{d_s}.
\end{equation}
Since $e^{-s|x|}\le1$, \eqref{eq:fs-detailed} yields $0\le f_s(z)\le sv(z)$.
If $Z$ is uniform on $Q_{d_s}$, then $T(Z)$ has the
$\operatorname{Gamma}(d_s,1)$ distribution. This allows us to control
the large values of $sv$ through Gamma tail probabilities.

We now consider the normalized first chaos with kernel $sv$ instead of $f_s$. Put
$\overline\sigma_s^{\,2}=s^3 3^{-d_s}$ and
$\overline X_s=I_1(sv)/\overline\sigma_s$. Since
$\int_{Q_{d_s}}s^2v(z)^2\dd z=s^23^{-d_s}$, the isometry
\eqref{eq:VarI1} shows that $\Var\overline X_s=1$. Moreover, by
\eqref{eq:V1s}, $\sigma_s^2=V_{1,s}\sim\overline\sigma_s^{\,2}$ for
$\beta>12/7$.

\begin{lemma}
	\label{lem:variance-kernel-replacement}
	If $d_s/\log s\to\beta>12/7$, then
	$
	\|X_s-\overline X_s\|_2\to0.
	$
\end{lemma}
\begin{proof}
	Since $0\le f_s\le sv$, we have $(sv-f_s)^2\le s^2v^2-f_s^2$.
	The isometry then yields
	\begin{displaymath}
		\E\bigl(I_1(sv)-I_1(f_s)\bigr)^2
		=s\int_{Q_{d_s}}(sv-f_s)^2\dd z
		\le s\int_{Q_{d_s}}(s^2v^2-f_s^2)\dd z
		=\overline\sigma_s^{\,2}-\sigma_s^{2}
		=o(\sigma_s^{2}).
	\end{displaymath}
	Hence, writing
	\begin{displaymath}
		X_s-\overline X_s
		=\frac{I_1(f_s)-I_1(sv)}{\sigma_s}
		+\left(\frac{\overline\sigma_s}{\sigma_s}-1\right)
		\frac{I_1(sv)}{\overline\sigma_s},
	\end{displaymath}
	we see that both terms converge to zero in $L^2$.
\end{proof}

For Borel $B\subset(0,\infty)$, define
\begin{displaymath}
	\overline\nu_s(B)=s\int_{Q_{d_s}}
	\one\{sv(z)/\overline\sigma_s\in B\}\dd z.
\end{displaymath}
This measure has unit second moment, since
$\int_0^\infty u^2\overline\nu_s(\!\dd u)
=s\int_{Q_{d_s}}s^2v(z)^2\dd z/\overline\sigma_s^{\,2}=1$.
Write
\begin{displaymath}
	\overline L_s(\eps)=\int_{(\eps,\infty)}
	u^2\overline\nu_s(\!\dd u).
\end{displaymath}
Recall that $\vartheta_s$ is defined in \eqref{eq:vartheta-s}
and $\Phi$ denotes the standard normal distribution function.

\begin{proposition}
	\label{prop:lindeberg}
	Suppose that $d_s/\log s\to\beta>12/7$ and
	$\vartheta_s\to\vartheta\in[-\infty,+\infty]$.
	Then $X_s\wto N(0,\Phi(-\vartheta))$.
	In particular, $X_s\wto N(0,1)$ if $\beta\in(12/7,\beta_L)$ and
	$X_s\wto0$ if $\beta>\beta_L$.
\end{proposition}

\begin{proof}
	By the definitions of $\overline\nu_s$ and $\overline\sigma_s$,
	\begin{align*}
		\overline L_s(\eps)
		&=s\int_{Q_{d_s}}\frac{s^2v(z)^2}{\overline\sigma_s^{\,2}}
		\one\left\{\frac{sv(z)}{\overline\sigma_s}>\eps\right\}\dd z\\
		&=\int_{Q_{d_s}}3^{d_s}v(z)^2
		\one\left\{T(z)<\log\frac{s}{\eps\overline\sigma_s}\right\}\dd z
		=\Prob{\operatorname{Gamma}(d_s,3)<\log\frac{s}{\eps\overline\sigma_s}}.
	\end{align*}
	Indeed, $3^{d_s}v(z)^2$ is a probability density on $Q_{d_s}$ under
	which $-\log(1-Z_i)$ are independent
	$\operatorname{Exp}(3)$ random variables, and hence
	$T(Z)\sim\operatorname{Gamma}(d_s,3)$. Since $\log\overline\sigma_s=\frac32\log s-\frac{d_s}2\log3$ and
	$1/\beta_L=\log3-2/3$, standardizing the threshold by the Gamma mean
	$d_s/3$ and standard deviation $\sqrt{d_s}/3$ gives
	\begin{displaymath}
		\frac{\log(s/(\eps\overline\sigma_s))-d_s/3}{\sqrt{d_s}/3}
		=\frac{\frac{d_s}{2}\left(\log3-\frac23\right)-\frac12\log s-\log\eps}
		{\sqrt{d_s}/3}
		=\vartheta_s-\frac{3\log\eps}{\sqrt{d_s}}.
	\end{displaymath}
	Hence the central limit theorem gives
	$\overline L_s(\eps)\to\Phi(\vartheta)$.
	
	Let $p_{d,r}(t)=\frac{r^d}{\Gamma(d)}\,t^{d-1}e^{-rt}$, $t>0$, be the
	$\operatorname{Gamma}(d,r)$ density. One easily checks that
	$p_{d,1}(t)=3^{-d}e^{2t}p_{d,3}(t)$ and
	$\sup_{t>0}p_{d,3}(t)\le3/\sqrt d$.
	Put $L=\log(s/(\eps\overline\sigma_s))$. If $L\le0$, then
	$sv/\overline\sigma_s\le\eps$, and hence
	$\int_{(\eps,\infty)}u\,\overline\nu_s(\!\dd u)=0$. If $L>0$, then
	\begin{align*}
		\int_{(\eps,\infty)}u\,\overline\nu_s(\!\dd u)
		&=\frac{s^2}{\overline\sigma_s}
		\E\left[v(Z)\one\{T(Z)<L\}\right]\\
		&=\frac{s^2}{\overline\sigma_s}3^{-d_s}
		\int_0^L e^t p_{d_s,3}(t)\dd t
		\le\frac{3}{\sqrt{d_s}}\frac{s^2}{\overline\sigma_s}
		3^{-d_s}e^L
		=\frac{3}{\eps\sqrt{d_s}}\to0.
	\end{align*}
	Thus, for any $L$,
	\begin{equation}
		\label{eq:large-jump-first-moment}
		\int_{(\eps,\infty)}u\,\overline\nu_s(\dd u)\to0.
	\end{equation}
	By the characteristic function formula for compensated Poisson integrals,
	\[
	\log\E e^{it\overline X_s}
	=\int_{(0,\eps]}(e^{itu}-1-itu)\,\overline\nu_s(\!\dd u)
	+\int_{(\eps,\infty)}(e^{itu}-1-itu)\,\overline\nu_s(\!\dd u).
	\]
	For $u\le\eps$, Taylor's formula gives
	$|e^{itu}-1-itu+t^2u^2/2|\le |t|^3\eps u^2/6$. Hence
	\[
	\int_{(0,\eps]}(e^{itu}-1-itu)\,\overline\nu_s(\!\dd u)
	=-\frac{t^2}{2}\bigl(1-\overline L_s(\eps)\bigr)+R_{s,\eps},
	\qquad
	|R_{s,\eps}|\le\frac{|t|^3\eps}{6},
	\]
	since $\int_{(0,\eps]}u^2\overline\nu_s(\!\dd u)=1-\overline L_s(\eps)\le1$.
	For $u>\eps$, we use $|e^{itu}-1-itu|\le2|t|u$. Therefore
	\[
	\Bigl|\log\E e^{it\overline X_s}+\frac{t^2}{2}\Phi(-\vartheta)\Bigr|
	\le\frac{t^2}{2}\bigl|\Phi(-\vartheta)-1+\overline L_s(\eps)\bigr|
	+\frac{|t|^3\eps}{6}
	+2|t|\int_{(\eps,\infty)}u\,\overline\nu_s(\!\dd u).
	\]
	For fixed $\eps>0$, the first term tends to zero since
	$\overline L_s(\eps)\to\Phi(\vartheta)$, while the third tends to zero
	by \eqref{eq:large-jump-first-moment}. Thus, letting first $s\to\infty$
	and then $\eps\downarrow0$ gives
	$\overline X_s\wto N(0,\Phi(-\vartheta))$. The same convergence for
	$X_s$ follows from Lemma~\ref{lem:variance-kernel-replacement}.
\end{proof}

\begin{proof}[Proof of Theorems~\ref{thm:main-gaussian}
	and \ref{thm:limtbetaL}]
	In the ranges $[0,\beta_0]$, $(\beta_0,12/7]$,
	and $(12/7,\beta_L]$, the results follow,
	respectively, from Propositions~\ref{prop:low-gaussian},
	\ref{prop:positive-gaussian}, and \ref{prop:lindeberg},
	the last along with Corollary~\ref{cor:transfer}.
\end{proof}

\section{The stable regime $\beta_L<\beta<\beta_c$}
\label{sec:above-lind-trans}

In this section, we prove Theorem~\ref{thm:stable-full}.
By Proposition~\ref{prop:lindeberg}, the normalized
compensator first chaos $X_s=I_1(f_s)/\sigma_s$ converges to zero
when $\beta>\beta_L$. A nondegenerate limit therefore requires
a smaller scale, determined by the large values of $f_s$.

Recall from \eqref{eq:v-T-def} that $v(z)=\prod_{i=1}^{d_s}(1-z_i)$,
$T(z)=-\log v(z)$, and $f_s(z)\le sv(z)$.
If $Z$ is uniform on $Q_{d_s}$, then
$T(Z)\sim\operatorname{Gamma}(d_s,1)$.
Following the kernel-replacement approach of
Subsection~\ref{sec:zero-saddle}, we first establish the stable
limit for the simpler kernel $sv$, then transfer it to $f_s$,
and finally show that the remaining terms $M_s$ and $U_s$
in \eqref{eq:N-EN-expression} are negligible.

\subsection{Limit behaviour of $I_1(sv)$}
\label{sec:gamma-normalization}

We first recall the function $\sI$ and the parameters
$r_\beta$ and $\alpha_\beta$ defined in
Section~\ref{sec:model-main-stat}.
Let $\sI(u)=u-1-\log u$, $u>0$,
be the Cram\'er rate function of the $\operatorname{Exp}(1)$
distribution. For $\beta>0$, let $r_\beta\in(0,1)$
be the unique solution of $\sI(r_\beta)=1/\beta$, and set
$\alpha_\beta=1/r_\beta-1$.
Since $\sI$ is strictly decreasing on $(0,1)$,
$r_\beta$ is continuous and strictly increasing in $\beta$.
Moreover, $\sI(1/3)=1/\beta_L$ and
$\sI(1/2)<1/\beta_{\mathrm{sp}}$.
Consequently, for $\beta\in(\beta_L,\beta_{\mathrm{sp}})$,
we have $r_\beta\in(1/3,1/2)$ and
$\alpha_\beta\in(1,2)$.
The identity $\sI(\log2/2)=1/\beta_c$ also gives
$r_{\beta_c}=\log2/2$. 
For $s>1$, we also define $r^\star_s\in(0,1)$
and $\alpha^\star_s>0$ by
$\sI(r^\star_s)=\log s/d_s$ and
$\alpha^\star_s=1/r^\star_s-1$, respectively.
Thus $r^\star_s\to r_\beta$ and
$\alpha^\star_s\to\alpha_\beta$ as $s\to\infty$,
whenever $d_s/\log s\to\beta>0$.

Let $\Gamma_{d_s}\sim\operatorname{Gamma}(d_s,1)$
with distribution function $F_{d_s}$ and density $g_{d_s}$.
For $s>1$, let $t_s^\star$ be the unique solution of
$F_{d_s}(t_s^\star)=1/s$.

\begin{lemma}
	\label{lem:stable-ts}
	If $d_s/\log s\to\beta\in(\beta_L,\beta_{\mathrm{sp}})$,
	then $t_s^\star/d_s\to r_\beta$. More precisely,
	\begin{equation}
		\label{eq:ts-second-order}
		t_s^\star=d_s r_s^\star
		+\frac{\log d_s}{2\alpha_s^\star}
		+\frac1{\alpha_s^\star}
		\log\bigl((1-r_s^\star)\sqrt{2\pi}\bigr)
		+o(1).
	\end{equation}
	Moreover, $sg_{d_s}(t_s^\star)\to\alpha_\beta$.
\end{lemma}

\begin{proof}	
	Let $t_s$, $s>1$, satisfy $t_s/d_s\to r\in(0,1)$.
	Putting $y=t_s-x$, we have
	\[\frac{F_{d_s}(t_s)}{g_{d_s}(t_s)}
	=\int_0^{t_s}\exp\Bigl\{(d_s-1)\log\Bigl(1-\frac y{t_s}\Bigr)+y\Bigr\}\dd y.\]
	By $\log(1-u)\le-u$, the integrand is bounded by
	$e^{-cy}$ for some $c>0$ and large $s$.
	Moreover, for each fixed $y$, the exponent converges to
	$-(r^{-1}-1)y$. Thus, by dominated convergence,
	\begin{equation}
		\label{eq:gamma-boundary-uniform}
		F_{d_s}(t_s)
		=\frac{g_{d_s}(t_s)}{r^{-1}-1}(1+o(1)).
	\end{equation}
	
	Since $r_s^\star\to r_\beta\in(0,1)$, applying
	\eqref{eq:gamma-boundary-uniform} and Stirling's formula
	at $t_s=d_s r_s^\star+z$ with $z=O(\log d_s)$ gives
	\begin{equation}
		\label{eq:gamma-local-log}
		F_{d_s}(d_s r_s^\star+z)
		=\frac{\exp\{-d_s\sI(r_s^\star)+\alpha_s^\star z\}}
		{(1-r_s^\star)\sqrt{2\pi d_s}}\,(1+o(1)).
	\end{equation}	
	By the definition of $r_s^\star$,
	$d_s\sI(r_s^\star)=\log s$.
	Taking $z=\pm C\log d_s$ with $C$ sufficiently large
	in \eqref{eq:gamma-local-log}, we obtain
	\[
	F_{d_s}(d_s r_s^\star-C\log d_s)
	<\frac1s<F_{d_s}(d_s r_s^\star+C\log d_s)
	\]
	for large $s$. Since $F_{d_s}$ is increasing and
	$F_{d_s}(t_s^\star)=1/s$, it follows that
	$t_s^\star-d_s r_s^\star=\sO(\log d_s)$.
	We may therefore apply \eqref{eq:gamma-local-log}
	with $z=t_s^\star-d_s r_s^\star$, obtaining
	\[
	\frac{\exp\{\alpha_s^\star(t_s^\star-d_s r_s^\star)\}}
	{(1-r_s^\star)\sqrt{2\pi d_s}}\,(1+o(1))=1.
	\]
	Taking logarithms proves \eqref{eq:ts-second-order}.
	In particular, $t_s^\star/d_s\to r_\beta$. Finally, the last claim now follows from \eqref{eq:gamma-boundary-uniform} at $t_s=t_s^\star$.
\end{proof}

Recall that $b_s=s e^{-t_s^\star}$. Since
\[
s\Prob{sv(Z)>b_s}
=s\Prob{T(Z)<t_s^\star}
=sF_{d_s}(t_s^\star)=1,
\]
where $Z$ is uniform on $Q_{d_s}$, this is the natural
normalization for the extreme values of $sv$.
Consider the normalized compensated Poisson integral
$\widetilde X_s=b_s^{-1}I_1(sv)$.
By the characteristic function formula for compensated
Poisson integrals,
\begin{equation}
\label{eq:chf}
\log\E e^{it\widetilde X_s}
=\int_0^\infty(e^{itu}-1-itu)\,\widetilde\nu_s(\!\dd u),
\end{equation}
where
\begin{equation}
\label{eq:hnu}
\widetilde\nu_s(B)
=s\int_{Q_{d_s}}\one\Bigl\{\frac{sv(z)}{b_s}\in B\Bigr\}\dd z=
s\int_{Q_{d_s}}\one\{e^{t_s^\star-T(z)}\in B\}\dd z,
\qquad B\in\mathcal B((0,\infty)).
\end{equation}

\begin{lemma}
	\label{lem:surrogate-stable}
	If $d_s/\log s\to\beta\in(\beta_L,\beta_{\mathrm{sp}})$,
	then $\widetilde X_s\wto Z_{\alpha_\beta}$,
	where $Z_{\alpha_\beta}$ is as in
	Theorem~\ref{thm:stable-full}.
\end{lemma}

\begin{proof}
	By \eqref{eq:hnu}, $\widetilde\nu_s$ is
	$s$ times the law of $e^{t_s^\star-\Gamma_{d_s}}$,
	where $\Gamma_{d_s}\sim\operatorname{Gamma}(d_s,1)$.
	Hence,
	\begin{equation*}
		\widetilde\nu_s(\!\dd u)
		=sg_{d_s}(t_s^\star)
		\left(1-\frac{\log u}{t_s^\star}\right)^{d_s-1}
		\one_{\{0<u<e^{t_s^\star}\}}\dd u.
	\end{equation*}
	
	By Lemma~\ref{lem:stable-ts},
	$sg_{d_s}(t_s^\star)\to\alpha_\beta$ and
	$(d_s-1)/t_s^\star\to1+\alpha_\beta\in(2,3)$.
	Hence, the density of $\widetilde\nu_s$ converges pointwise to
	$\alpha_\beta u^{-1-\alpha_\beta}$ and, by
	$\log(1-x)\le-x$, is bounded by
	$2\alpha_\beta u^{-(d_s-1)/t_s^\star}$
	for all large $s$.
	Together with
	$|e^{itu}-1-itu|\le\min\{t^2u^2/2,2|t|u\}$,
	this yields, for some $\varepsilon>0$ and all large $s$,
	\[
	|e^{itu}-1-itu|\frac{\dd\widetilde\nu_s}{\dd u}(u)
	\le C_t\bigl(u^{-1+\varepsilon}\one_{\{u\le1\}}
	+u^{-1-\varepsilon}\one_{\{u>1\}}\bigr),
	\]
	which is integrable on $(0,\infty)$. Dominated convergence therefore gives
	\[
	\log\E e^{it\widetilde X_s}
	\to\int_0^\infty
	(e^{itu}-1-itu)\alpha_\beta u^{-1-\alpha_\beta}\dd u.
	\]
	The limit is the characteristic exponent of $Z_{\alpha_\beta}$.	
\end{proof}

\subsection{Limit behaviour of $I_1(f_s)$}
\label{sec:uniform-replacement}

Now extend
$b_s^{-1}I_1(sv)\wto Z_{\alpha_\beta}$ to $b_s^{-1}I_1(f_s)\wto Z_{\alpha_\beta}$.
On $\{v>0\}$, let $r_s(z)=f_s(z)/sv(z)\le1$.

\begin{lemma}
	\label{lem:exact-first-chaos-stable}
	If $d_s/\log s\to\beta\in(\beta_L,\beta_{\mathrm{sp}})$,
	then $b_s^{-1}I_1(f_s)\wto Z_{\alpha_\beta}$.
\end{lemma}

\begin{proof}	
By \eqref{eq:fs-detailed}, $r_s(z)=\E e^{-s|X|}$,
where $X$ is uniform on $[z,1]$. Hence,
\begin{equation}
	\label{eq:rs}
	1-r_s(z)\le s\E|X|
	=s\prod_{i=1}^{d_s}\frac{1+z_i}{2}.
\end{equation}
Recall the generalized H\"older inequality
\[
\sum_{j=1}^m\left(\prod_{i=1}^d a_{ij}\right)^{1/d}
\le\prod_{i=1}^d\left(\sum_{j=1}^m a_{ij}\right)^{1/d},
\qquad a_{ij}\ge0.
\]
Applying it with $m=2$, $d=d_s$,
$a_{i1}=(1+z_i)/2$, and $a_{i2}=(1-z_i)/2$, we obtain
\[
\left(\prod_{i=1}^{d_s}\frac{1+z_i}{2}\right)^{1/d_s}
+\left(\prod_{i=1}^{d_s}\frac{1-z_i}{2}\right)^{1/d_s}
\le1.
\]
Since $\prod_{i=1}^{d_s}(1-z_i)=e^{-T(z)}$,
it follows from \eqref{eq:rs} that
\[
1-r_s(z)
\le s\left(1-\frac12e^{-T(z)/d_s}\right)^{d_s}.
\]
Since, by Lemma~\ref{lem:stable-ts}, $t_s^\star/d_s\to r_\beta$, we have for every fixed $C\in\R$
\[
\sup_{\{T(z)<t_s^\star+C\}}(1-r_s(z))
\le s\left(1-\frac12e^{-(t_s^\star+C)/d_s}\right)^{d_s}
=s^{1+\beta\log(1-e^{-r_\beta}/2)+o(1)}.
\]
The exponent on the right-hand side is negative, since
$1=\beta\sI(r_\beta)$ and
$\sI(u)+\log(1-e^{-u}/2)<0$ for $u\in(1/3,1/2)$.
Thus, for every $\eps>0$,
$\sup_{\{e^{t_s^\star-T(z)}>\eps\}}(1-r_s(z))\to0$.

By \eqref{eq:chf} and the analogous formula for
$I_1(f_s)$, using $f_s=r_s sv$, we have
\begin{equation*}
	\log\E e^{itI_1(f_s)/b_s}
	-\log\E e^{itI_1(sv)/b_s}
	=s\int_{Q_{d_s}}
	\left(e^{itr_s(z)sv(z)/b_s}
	-e^{itsv(z)/b_s}
	-it(r_s(z)-1)\frac{sv(z)}{b_s}\right)\dd z.
\end{equation*}
Since
$|e^{ita}-e^{itb}-it(a-b)|$ is bounded by both
$t^2|a^2-b^2|/2$ and $2|t||a-b|$, we have
\begin{align*}
	\bigl|\log\,&\E e^{itI_1(f_s)/b_s}
	-\log\E e^{itI_1(sv)/b_s}\bigr|\\
	&\le t^2s\int_{\{sv/b_s\le\eps\}}
	\biggl(\frac{sv(z)}{b_s}\biggr)^2\dd z
	+2|t|\sup_{\{sv/b_s>\eps\}}(1-r_s(z))
	\,s\int_{\{sv/b_s>\eps\}}
	\frac{sv(z)}{b_s}\dd z\\
	&=t^2\int_0^\eps u^2\widetilde\nu_s(\!\dd u)
	+2|t|\sup_{\{e^{t_s^\star-T(z)}>\eps\}}(1-r_s(z))
	\int_\eps^\infty u\,\widetilde\nu_s(\!\dd u).
\end{align*}

As shown in the proof of Lemma~\ref{lem:surrogate-stable},
the density of $\widetilde\nu_s$ is bounded by
$2\alpha_\beta u^{-(d_s-1)/t_s^\star}$
for large $s$.
Since $(d_s-1)/t_s^\star\to1+\alpha_\beta\in(2,3)$
by Lemma~\ref{lem:stable-ts}, the first integral
on the right-hand side tends to zero as $\eps\downarrow0$,
while the second remains bounded for each $\eps>0$,
both uniformly in large $s$.
Hence, letting first $s\to\infty$ and then $\eps\downarrow0$,
we conclude that
\[
\log\E e^{itI_1(f_s)/b_s}
-\log\E e^{itI_1(sv)/b_s}\to0,
\]
which extends $b_s^{-1}I_1(sv)\wto Z_{\alpha_\beta}$ from
Lemma~\ref{lem:surrogate-stable} to
$b_s^{-1}I_1(f_s)\wto Z_{\alpha_\beta}$.
\end{proof}

\subsection{Proof of Theorem~\ref{thm:stable-full}}
\label{sec:stable-remainder}

It remains to extend $b_s^{-1}I_1(f_s)\wto Z_{\alpha_\beta}$
to $b_s^{-1}(N_s-\E N_s)\wto Z_{\alpha_\beta}$.

\begin{lemma}
	\label{lem:stable-remainder}
	If $d_s/\log s\to\beta\in(\beta_L,\beta_{\mathrm{sp}})$,
	then $R_s/b_s^2\to0$.
	Moreover, $\E N_s/b_s^2\to0$ or $\infty$
	if $\beta<\beta_c$ or $\beta>\beta_c$, respectively.
\end{lemma}

\begin{proof}
	Lemmas~\ref{lem:mean-strong}, \ref{lem:R-bounds}, and
	\ref{lem:stable-ts} respectively yield
	$\E N_s=s^{2-\beta\log2+o(1)}$,
	$R_s\le s^{4-\beta\log6+o(1)}$, and
	$b_s^2=s^{2-2\beta r_\beta+o(1)}$.
	Hence,
	\[
	\frac{R_s}{b_s^2}
	\le s^{\beta(2\sI(r_\beta)+2r_\beta-\log6)+o(1)},
	\qquad
	\frac{\E N_s}{b_s^2}
	=s^{\beta(2r_\beta-\log2)+o(1)},
	\]
	where we used $1=\beta\sI(r_\beta)$.
	Since $r_\beta\in(1/3,1/2)$ and the function
	$u\mapsto2\sI(u)+2u-\log6$ is decreasing
	on this interval and equals
	$\log(3/2)-2/3<0$ at $u=1/3$,
	we have $R_s/b_s^2\to0$.
	The remaining limits follow from the strict
	monotonicity of $r_\beta$ and $r_{\beta_c}=\log2/2$.
\end{proof}

\begin{proof}[Proof of Theorem~\ref{thm:stable-full}]
	By Lemma~\ref{lem:exact-first-chaos-stable},
	$b_s^{-1}I_1(f_s)\wto Z_{\alpha_\beta}$
	for $\beta\in(\beta_L,\beta_{\mathrm{sp}})$.
	For $\beta\in(\beta_L,\beta_c)$,
	\eqref{eq:bound-difference-chaos} and
	Lemma~\ref{lem:stable-remainder} imply
	$b_s^{-1}(N_s-\E N_s-I_1(f_s))\to0$ in $L^2$.
	This completes the proof.
\end{proof}

\section{The Gaussian regime
$\beta_c<\beta<\beta_{\mathrm{sp}}$ and at $\beta_{\mathrm{sp}}$}
\label{sec:above-betac}

In the preceding regimes, the first chaos determined the limiting
distribution. Above $\beta_c$, this is no longer the case, and
we use martingale methods instead.

For each $s$, the processes $M_s(t)$ and $A_s(t)$,
$t\in[0,1]$, defined in \eqref{eq:MsAs}, are respectively
a square-integrable martingale w.r.t.\ the filtration
in \eqref{eq:Fst-def} and its predictable quadratic variation.
Moreover, $M_s(1)=M_s$, $A_s(1)=A_s$, and
$\E M_s^2=\E A_s=\E N_s$ by \eqref{eq:EM2}.

Our main tool is the following joint limit result
for predictable and deterministic Poisson integrals,
which will also be used in Section~\ref{sec:bound-behav-at-1}.
Set $\phi(t)=e^{it}-1-it$.

\begin{lemma}
	\label{lem:joint-clt}
	For each $s$, let $\xi_s=(\xi_s(x),x\in Q_{d_s})$ be bounded and predictable
	w.r.t.\	$(\sF_{s,t})_{t\in[0,1]}$, and let
	$j_s:Q_{d_s}\to[0,\infty)$ be deterministic with
	$j_s\in L^2(\mu_s)$. Suppose that
	\begin{enumerate}[label=\textup{(\alph*)}]
		\item $\|\xi_s\|_\infty\pto0$ and
		$s\int_{Q_{d_s}}\xi_s(x)^2\dd x\pto1$\rm;
		\item $s\int_{Q_{d_s}}\phi(vj_s(x))\dd x
		\to\Psi(v)$ for every $v\in\R$, where
		$e^\Psi$ is the characteristic function
		of a random variable $Z$;
		\item $s\int_{Q_{d_s}}|\xi_s(x)|j_s(x)\dd x\pto0$.
	\end{enumerate}
	Then
	\[
	\biggl(\int_{Q_{d_s}}\xi_s(x)\,\widehat\eta_s(\!\dd x),
	\int_{Q_{d_s}}j_s(x)\,\widehat\eta_s(\!\dd x)\biggr)
	\wto(G,Z),
	\]
	where $G\sim N(0,1)$ and $Z$ are independent.
\end{lemma}

\begin{proof}
	Fix $u,v\in\R$ and put $\chi_s=u\xi_s+vj_s$.
	For $t\in[0,1]$, let
	\[
	Z_s(t)=\int_{|x|\le t}\chi_s(x)\,
	\widehat\eta_s(\!\dd x),\qquad
	K_s(t)=s\int_{|x|\le t}
	\phi(\chi_s(x))\dd x.
	\]
	By the identity $\phi(a+b)=\phi(a)+\phi(b)+(e^{ia}-1)(e^{ib}-1)$,	we have
	\[K_s(1)
		=s\int_{Q_{d_s}}\phi(u\xi_s(x))\dd x
		+s\int_{Q_{d_s}}\phi(vj_s(x))\dd x+s\int_{Q_{d_s}}
		\bigl(e^{iu\xi_s(x)}-1\bigr)\bigl(e^{ivj_s(x)}-1\bigr)\dd x.\]
	Since $|\phi(a)+a^2/2|\le |a|^3/6$,
	condition (a) gives
	$s\int_{Q_{d_s}}\phi(u\xi_s)\dd x
	\pto-u^2/2$.
	Moreover, since $|(e^{iu\xi_s}-1)(e^{ivj_s}-1)| \le |uv\xi_s|\,j_s$, we have
	by (b) and (c)
	\begin{equation}
		\label{eq:joint-exponent}
		K_s(1)\pto -u^2/2+\Psi(v)=k(u,v).
	\end{equation}
	
	Denote
	\[\mathcal E_s(t)=\exp\{iZ_s(t)-K_s(t)\}=
	\exp\Bigl(\int_{|x|\le t}
	i\chi_s(x)\,\eta_s(\!\dd x)-s\int_{|x|\le t}(e^{i\chi_s(x)}-1)\dd x\Bigr).\]
	Since the first integral has finitely many jumps and
	the integrand of the second is bounded, $\mathcal E_s$
	has finite variation a.s. Applying the change of variables formula for processes
	of finite variation (see, e.g., \cite[Theorem~II.31]{P05})
	to $iZ_s-K_s$ and combining contributions from jumps and the continuous part,
	we obtain
	\[
		\mathcal E_s(t)=1+\int_{|x|\le t}
		\mathcal E_s(|x|-)(e^{i\chi_s(x)}-1)\,
		\widehat\eta_s(\!\dd x).
	\]
	Since the integrand is predictable and bounded
	for each fixed $s$, $\mathcal E_s(\cdot)$ is a
	square-integrable martingale.
	
	Moreover, by definition of $\mathcal E_s$, we have
	\[
	|\mathcal E_s(t)|=e^{-\operatorname{Re}K_s(t)},
	\qquad
	-\operatorname{Re}K_s(t)
	=s\int_{|x|\le t}(1-\cos\chi_s(x))\dd x.
	\]
	Thus $-\operatorname{Re}K_s(t)$ is continuous and nondecreasing.
	Fix $R>-\operatorname{Re}k(u,v)$ and set
	\[
	\tau_s=\inf\{t\in[0,1]:
	-\operatorname{Re}K_s(t)\ge R\}\wedge1.
	\]
	Hence, $|\mathcal E_s(\tau_s)|\le e^R$ by continuity of $-\operatorname{Re}K_s$.
	Since $\tau_s$ is a bounded stopping time, we have
	$\E\mathcal E_s(\tau_s)=\E\mathcal E_s(0)=1$.
	
	By \eqref{eq:joint-exponent},
	\begin{equation*}
	\Prob{\tau_s<1}
	\le\Prob{-\operatorname{Re}K_s(1)\ge R}\to0,
	\end{equation*}
	and so $K_s(\tau_s)\pto k(u,v)$.
	As $|e^{K_s(\tau_s)}|\le1$, this also yields
	$e^{K_s(\tau_s)}\to e^{k(u,v)}$ in $L^1$.
	Hence,
	\begin{equation}
		\label{eq:conv1}
		\bigl|\E e^{iZ_s(\tau_s)}-e^{k(u,v)}\bigr|=
		\bigl|\E\,\mathcal E_s(\tau_s)e^{K_s(\tau_s)}-\E\,\mathcal E_s(\tau_s)e^{k(u,v)}\bigr|
		\le e^R\,\E\bigl|	e^{K_s(\tau_s)}-e^{k(u,v)}\bigr|\to0.
	\end{equation}
	Moreover,
	\begin{equation}
	\label{eq:conv2}
	\bigl|\E e^{iZ_s(1)}-\E e^{iZ_s(\tau_s)}\bigr|\le2\Prob{\tau_s<1}\to0.
	\end{equation}
	Hence the claim follows from \eqref{eq:conv1}, \eqref{eq:conv2}, and L\'evy's continuity theorem. 
\end{proof}

We can now prove Theorems~\ref{thm:upper-gaussian}
and \ref{thm:densesp} (i) together.

\begin{proof}[Proof of Theorems~\ref{thm:upper-gaussian}
	and \ref{thm:densesp} (i)]
	Recall that
	$A_s-\E N_s=I_1(f_s)+U_s$, where
	$\E U_s^2=R_s$.
	In both cases, Lemma~\ref{lem:mean-strong} gives
	$\E N_s\to\infty$.
	Moreover, by Lemma~\ref{lem:R-bounds},
	$R_s/\E N_s\to0$.	
	
	We next show that
	$I_1(f_s)/\sqrt{\E N_s}\pto0$.
	If $\beta\in(\beta_c,\beta_{\mathrm{sp}})$,
	then $I_1(f_s)/b_s$ is tight by
	Lemma~\ref{lem:exact-first-chaos-stable},
	while $b_s^2/\E N_s\to0$ by
	Lemma~\ref{lem:stable-remainder}.
	If instead $\beta=\beta_{\mathrm{sp}}$,
	then \eqref{eq:V1s} yields
	$V_{1,s}/\E N_s\to0$.
	Thus $I_1(f_s)/\sqrt{\E N_s}\pto0$
	in both cases. Consequently,
	$(A_s-\E N_s)/\sqrt{\E N_s}\pto0$,
	and hence $A_s/\E N_s\pto1$.
	
	Let $\xi_s(x)=H_{s,x}/\sqrt{\E N_s}$.
	Since $H_{s,x}\in\{0,1\}$ by \eqref{eq:Hsx},
	\[
	\|\xi_s\|_\infty
	\le\frac1{\sqrt{\E N_s}}\to0,
	\qquad
	s\int_{Q_{d_s}}\xi_s^2\dd x
	=\frac{A_s}{\E N_s}\pto1.
	\]
	Applying Lemma~\ref{lem:joint-clt} with
	$j_s\equiv0$, we obtain
	$M_s/\sqrt{\E N_s}\wto N(0,1)$.
	Since $N_s-\E N_s=M_s+(A_s-\E N_s)$,
	we conclude that
	\[
	\frac{N_s-\E N_s}{\sqrt{\E N_s}}
	\wto N(0,1).
	\]
	
	Finally, if $\beta\in(\beta_1,\beta_{\mathrm{sp}}]$,
	then \eqref{eq:V1s} and
	Lemma~\ref{lem:mean-strong} give
	$V_{1,s}/\E N_s\to0$.
	Together with $R_s/\E N_s\to0$,
	Proposition~\ref{prop:variance-identity}
	yields 
	\begin{equation}
		\label{eq:var-mean-upper}
		\Var N_s\sim\E N_s,
		\qquad \beta\in(\beta_1,\beta_{\mathrm{sp}}].
	\end{equation}
	For $\beta\in(\beta_c,\beta_1)$,
	by \eqref{eq:NV} we have
	$\E N_s=o(\Var N_s)$, and hence $N_s-\E N_s$ normalized by its standard deviation converges
	to zero in probability.
	This completes both proofs.
\end{proof}

\section{The stable--Gaussian transition at
	\texorpdfstring{$\beta_c$}{beta c}}
\label{sec:bound-behav-at-1}

Now assume that $d_s/\log s\to\beta_c$.
By \eqref{eq:V1s} and Lemmas~\ref{lem:mean-strong},
\ref{lem:R-bounds}, and \ref{lem:stable-ts}, we have
\begin{equation}
  \label{eq:basic-asymptotics}
  \begin{aligned}
    \sqrt{\E N_s}&=s^{1-\beta_c\log2/2+o(1)},
    & V_{1,s}&=s^{3-\beta_c\log3+o(1)},\\
    R_s&\le s^{4-\beta_c\log6+o(1)},
    & b_s&=s^{1-\beta_c\log2/2+o(1)}.
  \end{aligned}
\end{equation}
In particular, $\sqrt{\E N_s}$ and $b_s$ have the same logarithmic order. We first establish a joint Gaussian--stable limit. Recall that $\alpha_c=\alpha_{\beta_c}=1/r_{\beta_c}-1=2/\log2-1$.

\begin{lemma}
  \label{lem:joint-boundary}
  If $d_s/\log s\to\beta_c$, then
  \[
    \left(
      \frac{M_s}{\sqrt{\E N_s}},
      \frac{I_1(f_s)}{b_s}
    \right)
    \wto(G,Z_{\alpha_c}),
  \]
  where $G\sim N(0,1)$ and the centered spectrally positive $\alpha_c$-stable
  variable $Z_{\alpha_c}$ are independent.
\end{lemma}

\begin{proof}
  The result follows from Lemma~\ref{lem:joint-clt}
  applied with $\xi_s(x)=H_{s,x}/\sqrt{\E N_s}$
  and $j_s=f_s/b_s$, once its conditions (a)--(c)
  are verified.
	
  Since $1-\beta_c\log2/2>0$,
  \eqref{eq:basic-asymptotics} yields
  $\E N_s\to\infty$, and so $\|\xi_s\|_\infty\le(\E N_s)^{-1/2}\to0$.
  Moreover,
  \[
    \E\biggl(\frac{A_s}{\E N_s}-1\biggr)^2
    =\frac{V_{1,s}+R_s}{(\E N_s)^2}\to0.
  \]
  Since $H_{s,x}^2=H_{s,x}$, we have
  $s\int_{Q_{d_s}}\xi_s^2\dd x
  =\frac{A_s}{\E N_s}\pto1$,	which proves (a).	
  Condition (b) is just a reformulation of
  Lemma~\ref{lem:exact-first-chaos-stable}
  in terms of characteristic functions.
  
  Using
  $\E H_{s,x}\le\E\eta_s([0,x])=s|x|$,
  $f_s(x)\le sv(x)$, and
  $\int_{Q_{d_s}}|x|v(x)\dd x=6^{-d_s}$,
  we obtain
  \begin{equation*}
    \E\Bigl(s\int_{Q_{d_s}}
    |\xi_s(x)|j_s(x)\dd x\Bigr)
    \le\frac{s^36^{-d_s}}
    {\sqrt{\E N_s}\,b_s}
    =s^{1-\beta_c\log3+o(1)}
    \to0,
  \end{equation*}
  which proves (c).
\end{proof}

\begin{proof}[Proof of Theorem~\ref{thm:boundary-trichotomy}]
  Recall from Subsection~\ref{sec:gamma-normalization}
  that $r_\beta\in(0,1)$ is defined by
  $\sI(r_\beta)=1/\beta$, with
  $r_{\beta_c}=\log2/2$, and that
  $r_s^\star$ and $\alpha_s^\star$ are defined by
  $\sI(r_s^\star)=\log s/d_s$ and
  $\alpha_s^\star=1/r_s^\star-1$.
  Recall also that $b_s=se^{-t_s^\star}$. Finally, set
  $\Delta_s=d_s-\beta_c\log s$.
  	
  By Lemma~\ref{lem:mean-strong},
  $\sqrt{\E N_s}=se^{-d_sr_{\beta_c}+o(1)}$.
  Hence \eqref{eq:ts-second-order} yields
  \[
    \log\frac{\sqrt{\E N_s}}{b_s}
    =t_s^\star-d_sr_{\beta_c}+o(1)
    =d_s(r_s^\star-r_{\beta_c})
    +\frac{\log d_s}{2\alpha_s^\star}
    +\frac1{\alpha_s^\star}
    \log\bigl((1-r_s^\star)\sqrt{2\pi}\bigr)
    +o(1).
  \]
  Since $\sI(r_\beta)=1/\beta$,
  differentiation gives
  $r_\beta'=1/(\alpha_\beta\beta^2)$.
  Hence, by Taylor's formula,
  \[
    d_s(r_s^\star-r_{\beta_c})
    =\frac{\Delta_s}{\alpha_c\beta_c}
    +O\left(\frac{\Delta_s^2}{\log s}\right),
    \qquad
    \alpha_s^\star-\alpha_c
    =O\left(\frac{|\Delta_s|}{\log s}\right).
  \]
  Using $\log d_s=\log\log s+\log\beta_c+
  O(|\Delta_s|/\log s)$, we obtain
  \[
    \log\frac{\sqrt{\E N_s}}{b_s}
    =\frac{\Delta_s}{\alpha_c\beta_c}
    +\frac{\log\log s}{2\alpha_c}
    +\frac1{\alpha_c}
    \log\bigl((1-r_{\beta_c})\sqrt{2\pi\beta_c}\bigr)
    +O\left(
      \frac{\Delta_s^2+|\Delta_s|\log\log s}{\log s}
    \right)+o(1).
  \]
  Since $\Delta_s=o(\log s)$ and
  $d_s-d_s^\star=\Delta_s+O(\log\log s)$,
  this finally yields
  \begin{equation*}
    \log\frac{\sqrt{\E N_s}}{b_s}
    =\frac{d_s-d_s^\star}{\alpha_c\beta_c}
    +o(1+|d_s-d_s^\star|).
  \end{equation*}
  Hence, $\sqrt{\E N_s}/b_s$ tends to $0$ in (i),
  to $e^{x/(\alpha_c\beta_c)}$
  in (ii), and to $\infty$ in (iii).
  
  By \eqref{eq:basic-asymptotics}, we have
  \begin{equation}
    \label{eq:UbN}
    \frac{\E U_s^2}{b_s^2}
    +\frac{\E U_s^2}{\E N_s}
    =\frac{R_s}{b_s^2}+\frac{R_s}{\E N_s}
    \le s^{2-\beta_c\log3+o(1)}
    \to0.
  \end{equation}
  For (i) and (ii), the decomposition
  \[
    \frac{N_s-\E N_s}{b_s}
    =\frac{\sqrt{\E N_s}}{b_s}
    \frac{M_s}{\sqrt{\E N_s}}
    +\frac{I_1(f_s)}{b_s}
    +\frac{U_s}{b_s},
  \]
  together with Lemma~\ref{lem:joint-boundary}
  and \eqref{eq:UbN}, yields the stated limits. For (iii), we similarly use
  \[
    \frac{N_s-\E N_s}{\sqrt{\E N_s}}
    =\frac{M_s}{\sqrt{\E N_s}}
    +\frac{b_s}{\sqrt{\E N_s}}
    \frac{I_1(f_s)}{b_s}
    +\frac{U_s}{\sqrt{\E N_s}}.
  \]
  The result again follows from Lemma~\ref{lem:joint-boundary}
  and \eqref{eq:UbN}, since $I_1(f_s)/b_s$ is tight.
\end{proof}

\section{De-Poissonization and the binomial model}
\label{sec:de-poiss-binom}

In this section, we transfer the Poisson limits to a
fixed-size binomial sample via direct $L^1$ coupling
estimates. For small dimensions we couple the Pareto
count, while in the remaining regimes we couple the
non-Pareto count itself. The binomial lower Gaussian
normalization is $\sqrt{V_{1,n}}$ throughout
$0<\beta<\beta_L$. For $\beta<\beta_0$, this differs
from the Poisson normalization, which was driven by
fluctuations in the total number of points.

As before, let $d_n\in\mathbb N$ satisfy
\begin{equation}
	\label{eq:bin-dimension}
	\frac{d_n}{\log n}\to \beta\in(0,\infty).
\end{equation}
For each $n$, let $X_1,X_2,\ldots$ be independent uniform points in
$Q_{d_n}$, and let $N_m^{(n)}$ and $K_m^{(n)}=m-N_m^{(n)}$ count,
respectively, the non-Pareto and Pareto-minimal points among
$X_1,\ldots,X_m$. Set
\[
  N_n^{\rm bin}=N_n^{(n)},\qquad
  K_n^{\rm bin}=K_n^{(n)},\qquad
  \mu_n=\E N_n^{\rm bin}.
\]
Let $\Pi_n\sim\operatorname{Poisson}(n)$ be independent of the
sequence $X_1,X_2,\ldots$. Then the sample stopped at $\Pi_n$ is a homogeneous Poisson
process of intensity $n$ on $Q_{d_n}$, so that
\[
  N_{\Pi_n}^{(n)}\overset{\rm d}=N_n^{\rm P},\qquad
  K_{\Pi_n}^{(n)}\overset{\rm d}=K_n^{\rm P},
\]
where the superscript ${\rm P}$ denotes the corresponding Poisson
counts. 

We begin with the error estimates for this coupling. For $d\ge1$ and $t\ge0$, put
\[
  J_d(t)=\int_{Q_d}e^{-t|x|}\dd x.
\]

\begin{lemma}
  \label{lem:binomial-coupling}
  For every $n$ and $d_n\ge 1$,
  \begin{align}
    \label{eq:binomial-N-coupling}
    \E\bigl|N_{\Pi_n}^{(n)}-N_n^{(n)}\bigr|
    &\le2\sqrt n\min\{1,n2^{-d_n}\},\\
    \label{eq:binomial-K-coupling}
    \E\bigl|K_{\Pi_n}^{(n)}-K_n^{(n)}\bigr|
    &\le2\sqrt{2n}\,J_{d_n}(n/8).
  \end{align}
\end{lemma}
\begin{proof}
  We fix $n$, write $d=d_n$, and suppress the superscript $(n)$ throughout.
  Let $D_m$ count the Pareto points among $X_1,\ldots,X_m$ that
  become non-Pareto after inserting $X_{m+1}$. Then
  \begin{equation}\label{eq:insertionid}
    N_{m+1}-N_m=\one_{\{X_{m+1}\text{ is non-Pareto}\}}+D_m,
    \qquad
    K_{m+1}-K_m=\one_{\{X_{m+1}\text{ is Pareto}\}}-D_m.
  \end{equation}
  For $m\ge1$, by exchangeability and conditioning on $X_1=x$, using
  $mu(1-u)^{m-1}\le\min\{1,mu\}$ on $[0,1]$, we have
  \begin{equation}\label{eq:ED}
    \E D_m=m\int_{Q_d}|x|(1-|x|)^{m-1}\dd x
    \le\min\{1,m2^{-d}\}.
  \end{equation}
  Also, since $1-(1-u)^m\le\min\{1,mu\}$ on $[0,1]$,
  \begin{displaymath}
  	\Prob{X_{m+1}\text{ is non-Pareto}}
  	=\int_{Q_d}\bigl(1-(1-|x|)^m\bigr)\dd x
  	\le \min\{1,m2^{-d}\}.
  \end{displaymath}
  Thus by \eqref{eq:insertionid},
  $0\le\E(N_{m+1}-N_m)\le2\min\{1,m2^{-d}\}$ for $m \ge 0$.
  Since $N_m$ is nondecreasing, telescoping yields, for $k\ge0$,
  \[
    \E|N_k-N_n|\le2|k-n|,
    \qquad
    \E|N_k-N_n|\le2^{-d}|k-n|(k+n-1).
  \]
  Averaging over $k=\Pi_n$ and using $\E|\Pi_n-n|\le\sqrt n$ yields
  the first bound in \eqref{eq:binomial-N-coupling}. For the second,
  the Cauchy--Schwarz inequality gives
  \[
    \E\bigl[|\Pi_n-n|(\Pi_n+n-1)\bigr]
    \le\sqrt{n\{n+(2n-1)^2\}}\le2n^{3/2},
  \]
  proving \eqref{eq:binomial-N-coupling}.

  For the Pareto count, \eqref{eq:insertionid} and \eqref{eq:ED} give,
  with $u=|x|$,
  \[
    \E|K_{m+1}-K_m| \le \int_{Q_d} [(1-u)^{m} + mu(1-u)^{m-1}]  \dd x
    =\int_{Q_d}(1-u)^{m-1}[1+(m-1)u]\dd x,
    \qquad m\ge1.
  \]
  If $m\ge2$ and $v=(m-1)u$, the integrand is bounded by
  \[
    (1-u)^{m-1}[1+(m-1)u] \le (1+v)e^{-v}\le2e^{-v/2}\le2e^{-mu/4}.
  \]
  Checking $m=0,1$ separately therefore gives
  $\E|K_{m+1}-K_m|\le2J_d(m/4)$ for every $m\ge0$.
  As $J_d$ is nonincreasing, triangle inequality yields
  \[
    \E|K_k-K_n|\le2|k-n|J_d(\min\{k,n\}/4).
  \]
  Now we again average over $k=\Pi_n$. For $0\le u\le1$, the
  Cauchy--Schwarz inequality gives
  \[
    \E\bigl[|\Pi_n-n|e^{-u\min\{\Pi_n,n\}/4}\bigr]
    \le\sqrt n\bigl(\E e^{-u\min\{\Pi_n,n\}/2}\bigr)^{1/2}.
  \]
  Plugging in the Poisson moment generating function and using
  $e^{-u/2}\le1-u/4$, we obtain
  \[
    \E e^{-u\min\{\Pi_n,n\}/2} \le e^{-un/2}+\E e^{-u\Pi_n/2}
    = e^{-un/2}+\exp\{n(e^{-u/2}-1)\}
    \le2e^{-un/4}.
  \]
  Integrating over $x$ proves \eqref{eq:binomial-K-coupling}.
\end{proof}

Two consequences will be used below. First, if $\beta>2$, then
\eqref{eq:binomial-N-coupling} and Lemma~\ref{lem:mean-strong} give
\begin{equation}
  \label{eq:binomial-mean-strong}
  \mu_n\sim n^22^{-d_n}\sim\E N_n^{\rm P},
\end{equation}
because $|\mu_n-\E N_n^{\rm P}|\le2n^{3/2}2^{-d_n}$.
Second, \eqref{eq:general-exponent} in Proposition~\ref{prop:laplace}
with $(m,\rho)=(1,0)$ gives 
\begin{equation}
  \label{eq:void-integral-exponent}
  nJ_{d_n}(cn)=n^{\kappa(\beta)+o(1)},\qquad
  \kappa(\beta)=\beta(1-\log\beta)<1,
  \qquad 0<\beta<1,
\end{equation}
for every fixed $c>0$. Indeed, by \eqref{eq:M}, $M_{1,0}(t)=1/(1-t)$,
so the minimizer in \eqref{eq:xi-variational} is $t=1-\beta$ in this
case, giving $\xi_{1,0}(\beta)=\beta-1-\beta\log\beta$.

Below, $f_n$, $V_{1,n}$ and $b_n$ denote the Poisson quantities
defined in \eqref{eq:fs-detailed}, \eqref{eq:V1} and \eqref{eq:bs},
respectively, evaluated at $s=n$ and $d=d_n$. Note that for either
$F_m^{(n)}=N_m^{(n)}$ or $F_m^{(n)}=K_m^{(n)}$, and any $a_n>0$, the
triangle inequality gives
\begin{equation}
  \label{eq:binomial-centered-coupling}
  \E\left|
    \frac{F_{\Pi_n}^{(n)}-\E F_{\Pi_n}^{(n)}}{a_n} -\frac{F_n^{(n)}-\E
      F_n^{(n)}}{a_n} \right|
  \le\frac{2}{a_n}\E\bigl|F_{\Pi_n}^{(n)}-F_n^{(n)}\bigr|.
\end{equation}

\begin{theorem}
  \label{thm:binomial-main}
  Assume $d_n/\log n\to \beta\in(0,\infty)$.
  \begin{enumerate}[label=\textup{(\roman*)}]
    \item If $0<\beta<\beta_L$, then
    \[
      \frac{N_n^{\rm bin}-\mu_n}{\sqrt{V_{1,n}}}\wto N(0,1).
    \]
    \item If $\beta_L<\beta<\beta_c$, then
    \[
      \frac{N_n^{\rm bin}-\mu_n}{b_n}\wto Z_{\alpha_\beta},
    \]
    where the stable law is defined in Theorem~\ref{thm:stable-full}.
    \item If $\beta_c<\beta<\beta_{\mathrm{sp}}$, then
    \[
      \frac{N_n^{\rm bin}-\mu_n}{\sqrt{\mu_n}}\wto N(0,1).
    \]
    \item If $\beta>\beta_{\mathrm{sp}}$, then
    $\mu_n\sim n^22^{-d_n}\to0$ and $N_n^{\rm bin}\pto0$.
  \end{enumerate}
\end{theorem}
\begin{proof}
  We first treat $0<\beta\le\beta_0$. Lemma
  \ref{lem:r3} together with
  \eqref{eq:r3-def} and \eqref{eq:V1-power} imply
  \[
    \frac{n\int f_n(z)^3\dd z}{V_{1,n}^{3/2}}
    =n^{r_3(\beta)-3r(\beta)/2+o(1)}\to0.
  \]
  Thus, recalling $V_{1,n}$ from \eqref{eq:V1}, Lemma
  \ref{lem:poisson-lyapunov} yields
  \begin{equation}\label{eq:Iclt}
    \frac{I_1(f_n)}{\sqrt{V_{1,n}}}\wto N(0,1).
  \end{equation}
  For the Poisson Pareto count, write
  \[
    K_n^{\rm P}= \int_{Q_{d_n}}(1-H_{n,x})\,\eta_n(\!\dd x) =\overline M_n+n-A_n,\qquad
    \overline M_n=\int_{Q_{d_n}}(1-H_{n,x})\,
    \widehat\eta_n(\!\dd x).
  \]
  The predictable-integral isometry gives
  $\E\overline M_n^2=\E K_n^{\rm P}$. Since $\E U_n^2 = R_n$ by
  \eqref{eq:Us}, we obtain 
  \begin{equation}
    \label{eq:KsBin}
    \|K_n^{\rm P}-\E K_n^{\rm P}+I_1(f_n)\|_2
    \le \sqrt{\E \overline M_n^2} + \sqrt{\E U_n^2}= \sqrt{\E K_n^{\rm P}}+\sqrt{R_n}.
  \end{equation}
  Here $R_n=o(V_{1,n})$ by Lemma~\ref{lem:R-bounds}. To compare the
  exponents of $\E K_n^{\rm P}=nJ_{d_n}(n)=n^{\kappa(\beta)+o(1)}$
  with $V_{1,n}=n^{r(\beta)+o(1)}$ from \eqref{eq:V1-power}, denote
  $u_\beta:=3-2t_\beta>0$, where recall that $t_\beta \in [0,3/2)$
  given at \eqref{eq:tb} is the saddle of $r(\beta)$.
  Equations \eqref{eq:tb} and \eqref{eq:rbeta-saddle} give
  \[
    \beta=\frac{u_\beta(1+u_\beta)}{1+2u_\beta},\qquad
    r(\beta)-\kappa(\beta)
    =\frac{u_\beta^2}{1+2u_\beta}
     +\beta\log\frac4{1+2u_\beta}>0.
  \]
  The last inequality follows from $\beta\le\beta_0<3/8$, which
  implies $0<u_\beta<1/2$.
  Thus both terms on the right-hand side of
  \eqref{eq:KsBin} are $o(\sqrt{V_{1,n}})$, and
  $(K_n^{\rm P}-\E K_n^{\rm P})/\sqrt{V_{1,n}}\wto N(0,1)$ by \eqref{eq:Iclt}.
  Moreover, since $\kappa(\beta)<r(\beta)$
  and $\kappa(\beta)<1$, we have $2\kappa(\beta)-1<r(\beta)$. Thus,
  \eqref{eq:binomial-K-coupling},  
  \eqref{eq:void-integral-exponent} and \eqref{eq:V1-power} imply
  \[
    \frac{\E|K_{\Pi_n}^{(n)}-K_n^{(n)}|}{\sqrt{V_{1,n}}}
    \le n^{\kappa(\beta)-1/2-r(\beta)/2+o(1)}\to0.
  \]
  By \eqref{eq:binomial-centered-coupling}, the CLT for
  $K_n^{\rm bin}$, and therefore for $N_n^{\rm bin}=n-K_n^{\rm bin}$
  in (i) follows for $\beta \le \beta_0$.

  In the remaining nonsparse regimes, we couple $N_n^{(n)}$. The bound
  \eqref{eq:binomial-N-coupling} along with the asymptotics
  \eqref{eq:V1-power}, \eqref{eq:V1s}, Lemma~\ref{lem:stable-ts},
  and \eqref{eq:binomial-mean-strong} give
  $$
  \frac{\E|N_{\Pi_n}^{(n)}-N_n^{(n)}|}{a_n}
  \le n^{-\gamma(\beta)+o(1)},
  $$
  with
  \begin{equation}
    \label{eq:binomial-coupling-rates}
    (a_n,\gamma(\beta))=
    \begin{cases}
      \bigl(\sqrt{V_{1,n}},\,(r(\beta)-1)/2\bigr),
      &\beta_0<\beta\le12/7,\\[1mm]
      \bigl(\sqrt{V_{1,n}},\,\beta\log(4/3)/2\bigr),
      &12/7<\beta\le\beta_L,\\[1mm]
      \bigl(b_n,\,\beta(\log2-r_\beta)-1/2\bigr),
      &\beta_L<\beta\le\beta_c,\\[1mm]
      \bigl(\sqrt{\mu_n},\,(\beta\log2-1)/2\bigr),
      &\beta_c\le\beta\le\beta_{\mathrm{sp}}.
    \end{cases}
  \end{equation}
  Observe that every exponent $\gamma(\beta)$ is positive. For the
  first row in \eqref{eq:binomial-coupling-rates} this was shown in
  the proof of Proposition~\ref{prop:variance-dominance}; for the
  third, this follows upon noting that $r_\beta\le r_{\beta_c}=\log2/2$ and
  $\beta_L\log2>1$. The other rows are immediate.
  The conclusions of Theorems \ref{thm:main-gaussian}, \ref{thm:stable-full} and
  \ref{thm:upper-gaussian} can now be transferred via
  \eqref{eq:binomial-centered-coupling}. The assertions (i)--(iii) now
  follow immediately upon using $\Var N_n^{\rm P}\sim V_{1,n}$ from 
  Proposition~\ref{prop:variance-dominance} for (i), and
  \eqref{eq:binomial-mean-strong} for (iii).
  Finally, (iv) follows from \eqref{eq:binomial-mean-strong} and Markov's
  inequality.
\end{proof}

The same estimates also transfer the boundary results in Theorems
\ref{thm:limtbetaL}, \ref{thm:boundary-trichotomy} and
\ref{thm:densesp}.
Let $\vartheta_n$ and $d_n^\star$ be the quantities defined in
\eqref{eq:vartheta-s} and \eqref{eq:dsstar} in
Section~\ref{sec:model-main-stat}, with $s$ replaced by $n$.

\begin{corollary}
  \label{cor:binomial-critical}
  Assume $d_n/\log n\to\beta\in(0,\infty)$.
  \begin{enumerate}[label=\textup{(\roman*)}]
    \item If $\beta=\beta_L$ and
    $\vartheta_n\to\vartheta\in[-\infty,+\infty]$, then
    \[
      \frac{N_n^{\rm bin}-\mu_n}{\sqrt{V_{1,n}}}
      \wto N(0,\Phi(-\vartheta)).
    \]
    \item If $\beta=\beta_c$, all three conclusions of
    Theorem~\ref{thm:boundary-trichotomy} hold with
    $N_s-\E N_s$, $b_s$, $\sqrt{\E N_s}$ and $d_s-d_s^\star$
    replaced, respectively, by
    $N_n^{\rm bin}-\mu_n$, $b_n$, $\sqrt{\mu_n}$ and
    $d_n-d_n^\star$. In particular, if $d_{n_j}-d_{n_j}^\star\to x\in\R$
    along some sequence $n_j\to\infty$, then
    \[
    \frac{N_{n_j}^{\rm bin}-\mu_{n_j}}{b_{n_j}}
    \wto Z_{\alpha_c}+e^{x/(\alpha_c\beta_c)}G,
    \]
    with $G\sim N(0,1)$ independent of $Z_{\alpha_c}$.
    \item Let $\beta=\beta_{\mathrm{sp}}$ and
    $h_n=d_n-\beta_{\mathrm{sp}}\log n$. Then
    \begin{align*}
    	h_n\to-\infty
    	&\quad\Longrightarrow\quad
    	\frac{N_n^{\rm bin}-\mu_n}{\sqrt{\mu_n}}\wto N(0,1),\\
    	h_{n_j}\to x\in\R
    	&\quad\Longrightarrow\quad
    	N_{n_j}^{\rm bin}\wto\operatorname{Poisson}(2^{-x})
    	\quad\text{along }n_j\to\infty,\\
    	h_n\to+\infty
    	&\quad\Longrightarrow\quad N_n^{\rm bin}\pto0.
    \end{align*}
  \end{enumerate}
\end{corollary}
\begin{proof}
  At $\beta_L$, the result follows from Theorem~\ref{thm:limtbetaL}
  and \eqref{eq:binomial-centered-coupling}, using the second row in
  \eqref{eq:binomial-coupling-rates} above and
  $\Var N_n^{\rm P}\sim V_{1,n}$ from Proposition
  \ref{prop:variance-dominance}. At $\beta_c$, the third and fourth 
  rows of \eqref{eq:binomial-coupling-rates} transfer
  Theorem~\ref{thm:boundary-trichotomy} on both 
  normalizations $b_n$ and $\sqrt{\mu_n}$. At $\beta_{\mathrm{sp}}$,
  the fourth row of \eqref{eq:binomial-coupling-rates} applies 
  on the dense side $h_n \to -\infty$; here we use
  $\Var N_n^{\rm P}\sim\E N_n^{\rm P}\sim\mu_n$ by
  \eqref{eq:var-mean-upper}. For the remaining cases, 
  \eqref{eq:binomial-N-coupling} gives
  $$
  \E|N_{\Pi_n}^{(n)}-N_n^{(n)}| \le 2n^{3/2} 2^{-d_n}
  = 2n^{3/2} 2^{-\beta_{\mathrm{sp}} \log n -h_n} 
  =2n^{-1/2+o(1)}\to0.
  $$
  Hence the limits transfer from Theorem~\ref{thm:densesp} via
  Slutsky's theorem, completing the proof. 
\end{proof}

\begin{remark}
  The binomial Gaussian regime has the single scale
  $\sqrt{V_{1,n}}$ throughout $0<\beta<\beta_L$, including
  $\beta=\beta_0$. This is because fixing the sample size removes the
  Poisson total-count
  fluctuation. Also, the assumption $\beta>0$ is essential; for
  example, $d_n=1$ gives $N_n^{\rm bin}\equiv n-1$.
  The coupling estimates control centered counts in $L^1$ and do not
  by themselves establish binomial variance asymptotics. The scales
  in the theorem and corollary above should therefore not be replaced by
  $\sqrt{\Var N_n^{\rm bin}}$ without a separate variance argument.
\end{remark}

\section*{Acknowledgements}

CB was supported by the German Research Foundation (DFG) Project
531540467. During the preparation of this manuscript, the authors used
ChatGPT (OpenAI) to assist with language editing, improving the
clarity and conciseness of the exposition, and checking selected
calculations. All AI-assisted suggestions were reviewed by the
authors, who take full responsibility for the mathematical content and
the final manuscript.


\end{document}